\documentclass[a4paper,12pt,twoside]{amsart} 
\usepackage{amsmath,amsthm}
\usepackage{amsfonts}
\usepackage{amssymb}
\usepackage[english]{babel}

\newtheorem{theo}{Theorem}[section]
\newtheorem{lem}[theo]{Lemma}
\newtheorem{prop}[theo]{Proposition}
\newtheorem{cor}[theo]{Corollary}

\newcommand{\e}{{\text{e}}}

\begin{document}

\title[  Hyperbounded Markov operators]
{$L^2$-Quasi-compact  and  hyperbounded\\ Markov operators}

\author{Guy Cohen}
\address{School of Electrical Engineering, Ben-Gurion University, Beer-Sheva, Israel}
\email{guycohen@bgu.ac.il}

\author{Michael Lin}
\address{Department of Mathematics, Ben-Gurion University, Beer-Sheva, Israel}
\email{lin@math.bgu.ac.il}


\begin{abstract}
A Markov operator $P$ on a probability space $(S,\Sigma,\mu)$, with $\mu$ invariant,
is called {\it hyperbounded} if for some $1 \le p<q \le \infty$ it maps (continuously)
$L^p$ into $L^q$.

We deduce from a recent result of Gl\"uck that a hyperbounded $P$ is quasi-compact,
hence uniformly ergodic, in all $L^r(S,\mu)$, $1<r< \infty$. We prove, using a method similar 
to Foguel's, that a hyperbounded Markov operator has periodic behavior similar to that of 
Harris recurrent operators, and for the ergodic case obtain conditions for aperiodicity.

Given a probability $\nu$ on the unit circle, we prove  that if the convolution operator
$P_\nu f:=\nu*f$ is hyperbounded, then $\nu$ is atomless. We show that there is $\nu $
absolutely continuous such that $P_\nu$ is not hyperbounded, and there is  $\nu$ with all 
powers singular such that $P_\nu$ is hyperbounded. As an application, we prove that if 
$P_\nu$ is  hyperbounded, then for any sequence $(n_k)$ of distinct positive integers with
 bounded gaps, $(n_kx)$ is uniformly distributed mod 1 for $\nu$ almost every $x$ 
(even when $\nu$ is  singular).

\end{abstract}

\dedicatory{Dedicated to the memory of Shaul Foguel}
\maketitle

\section{Introduction}

Let $(S,\Sigma,\mu)$ be a probability space and $P(x,A)$ a transition probability on 
$S\times\Sigma$ with Markov operator $Pf(x):= \int f(y)P(x,dy)$ for  bounded measurable $f$. 
We call $P$ {\it bi-stochastic} on $(S,\Sigma,\mu)$ when  $\mu$ is {\it invariant} for $P$, i.e. 
$\int P(x,A)d\mu(x) = \mu(A)$ for every $A \in \Sigma$. In that case, if $f$ is bounded and 
$f=0$ a.e. $\mu$, then $Pf=0$ a.e., so $P$ defines an operator (still denoted by $P$) on
$L^\infty(\mu)$, and the invariance of $\mu$ yields
$$
\|Pf\|_1:= \int Pf(x) d\mu(x) = \int f d\mu = \|f\|_1  \quad \text{for } 0 \le f\in L^\infty(\mu).
$$
In this case $P$ clearly extends to a contraction of $L^1(\mu)$, and then we have that 
$P$ is a contraction of each $L^p(\mu)$, $1\le p \le \infty$ \cite[p. 65]{Kr} (or see 
\cite[Corollary VI.10.12]{DS}), and for $1\le p <\infty$ the power averages 
$\{A_n(P) := \frac1n\sum_{k=0}^{n-1} P^k \}_{n\ge 1}$ converge in the strong operator topology
of $L^p(\mu)$, and a.e.  \cite[Theorem 1.7.2]{Kr}.  We say that $P$ is {\it ergodic} if
 $Pf=f \in L^\infty(\mu)$ implies that $f$ is constant a.e. When $P$ is ergodic and bi-stochastic, 
$\lim A_n(P)f= \int f\,d\mu$ a.e. and in $L^p$-norm, for any $f \in L^p(\mu)$, $1 \le p < \infty$.
\smallskip

{\bf Definition.} Let $(S,\Sigma,\mu)$ be a probability space and $1 \le p < \infty$. 
A bounded operator $T$ on $L^p(S,\mu)$ is called {\it hyperbounded} if 
for some $q>p$ the operator $T$ maps $L^p(S,\mu)$ into $L^q(S,\mu)$.
As observed in \cite{G2}, a hyperbounded $T$ maps $L^p$ to $L^q$ continuously, 
by the closed graph theorem.
Note that if $T$ maps $L^p$ to $L^\infty$, then it maps $L^p$ to $L^q$ for any 
$p<q< \infty$, since $\|Tf\|_q \le \|Tf\|_\infty \le C\|f\|_p$.
\smallskip

Gl\"uck \cite[Theorem 1.1]{G2} proved the following.
\begin{theo} \label{gluck}
Let $1<p< \infty$ and let $T$ be a power-bounded positive operator on $L^p(S,\mu)$.
If $T$ is hyperbounded, then the essential spectral radius of $T$ (as an operator
on the complex $L^p$) is less than 1.
\end{theo}

A  Markov operator $P$ on a probability space $(S,\Sigma,\mu)$ is called {\it hyperbounded} 
if $\mu$ is invariant (i.e. $P$ is bi-stochastic), and for some $1< p <q \le \infty$ the operator 
$P$ maps $L^p(\mu)$ to $L^q(\mu)$ (i.e. $P$ is hyperbounded on $L^p(\mu)$). 
Since all $L^p(S,\mu)$ spaces  are invariant under $P$ and therefore under all its powers, 
it follows that {\it all the powers of a hyperbounded Markov operator are hyperbounded}.
\medskip

We deduce from Gl\"uck's result (Theorem \ref{gluck}) that a hyper-bounded bi-stochastic
Markov operator is uniformly ergodic in all $L^r(\mu)$ spaces, $1<r < \infty$. 
We prove, using a method similar to Foguel's \cite{F1}, \cite{F3}, that an ergodic 
hyperbounded Markov operator has a periodic behavior similar to that of ergodic Harris 
recurrent  operators, and obtain conditions for aperiodicity.
\smallskip 

A probability $\nu$ on the unit circle $\mathbb T$ defines a Markov operator $P_\nu$ by
$P_\nu f =\nu*f$, for which the normalized Lebesgue measure $\mu$ on $\mathbb T$ is invariant.
We prove that if $P_\nu$ is hyperbounded, then $\nu$ has no atoms. We show that there 
exists an absolutely continuous $\nu$ (so $P_\nu$ is uniformly ergodic in $L^2(\mathbb T)$ 
by \cite[Theorem 4.4]{DL}) which is not hyperbounded, and that there are singular 
$\nu$ such that $P_\nu$ is hyperbounded and not Harris recurrent.

\medskip

\section{Uniformly ergodic positive operators on reflexive Banach lattices}

A bounded linear operator $T$ on a Banach space $X$ is called {\it uniformly ergodic} if its 
power averages $\{A_n(T) \}_{n\ge 1}$ converge in the operator norm topology. An operator 
$T$ is uniformly ergodic if and only if $n^{-1}\|T^n\| \to 0$ and $(I-T)X$ is closed \cite{L1}.
When $X$ is over $\mathbb C$,  the powers $\{T^n\}$ converge in operator norm if and only 
 if $T$ is uniformly ergodic and $\sigma(T)\cap \mathbb T \subset \{1\}$ \cite[Theorem 4]{Lu};
in that case, $T$ is power-bounded and $r(T_{|(I-T)X}) <1$. Whether $X$ is over $\mathbb R$ or 
over $\mathbb C$, operator norm convergence of $T^n$, say to some $E$, is exponential: there 
exist $C>0$ and $\rho<1$ such that $\|T^n-E\| \le C\rho^n$ (see \cite[Proposition 3.1]{DL}).
\smallskip

A bounded linear $T$ on $X$ is called {\it quasi-compact} if there exists a compact operator $K$
such that $\|T^n-K\|<1$ for some $n \ge 1$. If $T$ is quasi-compact and $\frac1n T^n \to 0$
in the weak operator topology, then $T$ is uniformly ergodic, $\sigma(T) \cap \mathbb T$ is 
finite, and each $\lambda \in \sigma(T) \cap \mathbb T$ is a simple pole of the resolvent
(hence an eigenvalue) with finite-dimensional eigenspace  \cite[p. 711]{DS}. It is known 
\cite[Lemma 2.2.4, p. 88]{Kr} that $T$ is quasi-compact if and only if there is a sequence 
$\{K_n\}$ of compact operators such that $\|T^n -K_n\| \to 0$; consequently the powers of a 
quasi-compact operator are quasi-compact. Conversely, if $T^m$ is quasi-compact for some
 $m>1$, then $T$ is quasi-compact.

\begin{prop} \label{qc}
Let $T$ be a {\rm positive} power-bounded operator on a complex Banach lattice $L$.
If $T$ is quasi-compact, then there exists an integer $d \ge 1$ such that each  (of the 
finitely many) $\lambda \in \sigma(T) \cap \mathbb T$ is a {\rm d}th root of unity, 
$\sigma(T^d) \cap \mathbb T  \subset \{1\}$, and $(T^{nd})$ converges in operator
norm as $n \to \infty$. Moreover, each $\lambda \in \sigma(T) \cap \mathbb T$ is
an eigenvalue, with finite-dimensional eigenspace.
\end{prop}
\begin{proof}
Power-boundedness implies $r(T) \le 1$, so we have to prove only when $r(T)=1$.
\smallskip

By  \cite[Theorem VIII.8.3]{DS},  the peripheral spectrum $\sigma(T) \cap \mathbb T$
consists of finitely many points, which are simple poles, hence eigenvalues, 
with finite-dimensional corresponding eigenspaces. 
Since $T$ is power-bounded and positive, by a result of Lotz \cite[p. 327, Theorem 4.9]{Sc}
the peripheral spectrum  $\sigma(T) \cap \mathbb T$ is cyclic, i.e.\ $\lambda^n$ is in the
 peripheral spectrum when $\lambda$ is. 
Since $ \sigma(T)\cap \mathbb T$ is finite, $ \lambda \in \sigma(T) \cap \mathbb T$
must be a root of unity. Thus $\sigma(T)\cap \mathbb T$ is a finite set of roots of unity, 
and let $d$ be the smallest common multiple of their orders.

By the spectral mapping theorem, we have $\sigma(T^d) \cap \mathbb T  = \{1\}$.
Since $T^d$ is also quasi-compact, it is uniformly ergodic, so $(T^{nd})$
converges in operator norm, as $n \to \infty$, by \cite[Theorem 4]{Lu}.
\end{proof}

{\bf Remarks.} 1. For the  cases $L=C(S)$ for a compact Hausdorff space $S$ or
(by duality) $L=L^1(\mu)$, the existence of $d$ is proved in \cite[Lemma VIII.8.5]{DS}.

2. Without positivity the proposition is false, even in finite-dimensional spaces.


\begin{cor} \label{bartoszek}
Let $T$ be a positive power-bounded operator on a  complex Banach lattice $L$.
Then $T$ is quasi-compact if and only if for some  integer $d\ge 1$ and a finite-dimensional
projection $E$ we have $\|T^{nd}-E\| \to 0$.
\end{cor}

{\bf Remark.} For contractions the corollary was proved by Bartoszek \cite[Theorem 2]{Ba},
who gave a representation of the limit (see Section 4).
\smallskip

\begin{prop} \label{UE}
Let $T$ be a positive power-bounded operator on a  complex Banach lattice $L$.
If $T$ is uniformly ergodic and  $F(T):= \{ f \in L:\, Tf=f \}$ is finite-dimensional,
then $T$ is quasi-compact, and has all the properties of Proposition \ref{qc}.
\end{prop}
\begin{proof}
%
Since $T$ is positive uniformly ergodic with $F(T)$ is finite-dimensional,  
by \cite[Theorem 1]{L2} $T$ is quasi-compact. 
The existence of $d$ and the other properties follow from  Proposition \ref{qc}.
\end{proof}

{\bf Remark.} Without positivity the proposition is false (take\  $-I$ on $L$ infinite-dimensional); moreover, $T^2$ need not be uniformly ergodic when $T$ is \cite{L1a}.
\smallskip

{\bf Definition.}    Let  $1 \le p< \infty$.   A power-bounded   linear   operator $T$ on
$L^p(S,\mu)$ is said to be {\it uniformly integrable} in $L^p$ (in short, p-UI),
if $\{|Tf|^p:\, \|f\|_p \le 1\}$ is uniformly integrable.

When $p=1$, uniform integrability of $T$ means that $T$ is weakly compact; hence $T^2$
is compact in $L^1$ \cite[Corollary VI.8.13]{DS}, and  $F(T)\subset F(T^2)$ is
finite-dimensional. However, $T$ may be weakly compact in $L^1$ without being compact
\cite{YMK}.

Wu \cite[Proposition 1.4]{Wu} proved that if $T$ is compact or hyperbounded in $L^p$,
 then $T$ is p-UI. Gl\"uck \cite[Corollary 2.3]{G2} proved that if $T$ is a {\it positive} 
power-bounded operator on $L^p$ such that for some integer $m \ge 1$ $T^m$ is p-UI, 
then the essential spectral radius of $T$ is less than 1.
\smallskip

{\bf Remarks.} 1. An  example in \cite[Remark 1.1(2)]{BWY} shows that we may have a 
bi-stochastic Markov operator $P$ with $P^n$ compact in $L^2$ (so $P^n$ is 2-UI) for 
every large n, but $P $ itself is not 2-UI (take $P=P_{t_1}$ for some $t_1\in (0,r_0)$ there;
this $P$ is quasi-compact, hence uniformly ergodic, in $L^2$, but not 2-UI). 
Theorem 1.1 of \cite{BWY} yields that in that example $P_t$ is even hyperbounded in $L^2$ 
for large $t$, so $P$ as above is not hyperbounded in $L^2$, though $P^n$ is for large $n$.
By Proposition \ref{all-Lp} below, for any $1<p< \infty$, $P$ is not hyperbounded in $L^p$,
though $P^n$ is for large $n$.

2. El Machkouri, Jakubowski and Voln\'y \cite{MJV} have an example of a 2-UI bi-stochastic
Markov operator which is not hyperbounded in $L^2$.

\begin{prop} \label{fixed-space}
Let $1 < p< \infty$ and let $T$ be a positive power-bounded operator on $L^p(S,\Sigma,\mu)$.
Set $Ef = \lim_n \frac1n\sum_{k=1}^n T^kf$, and assume that $T$ is p-UI.

(i) If $Tf \ge f \ge 0$ implies $Tf=f$ a.e., then $F(T)$ is finite-dimensional.

(ii) If $Ef \not\equiv 0$ for  $0 \le f \not\equiv 0$, then $F(T)$ is finite-dimensional.

(iii) If $\|T\| \le 1$, then $F(T)$ is finite-dimensional.
\end{prop}
\begin{proof} If $E=0$ we have $F(T)=\{0\}$; we therefore assume $E \ne 0$.

 (i) If $g \in F(T)$, then $|g|=|Tg| \le T|g|$, hence $|g| \in F(T)$, so
$F(T)$ is a sublattice of $L^p(\mu)$. Then by \cite[Corollary 4.3]{BL}, $F(T)$ is the range
of a positive contractive projection. Hence by \cite[Theorem 4.1]{BL},  there is a 
measure space $(S',\Sigma',\nu)$ such that
$F(T)$ is isometrically isomorphic to $L^p(S', \Sigma',\nu)$. Since $T$ is p-UI and on $F(T)$ 
it is the identity, we obtain that the identity on $L^p(\nu)$ is p-UI, which means that 
the unit ball of $L^1(\nu)$ is weakly compact. This means that $L^1(\nu)$ is reflexive, 
so it is finite-dimensional. Hence $F(T)$ is finite-dimensional.

(ii) The assumption implies (i), by  the proof of \cite[Proposition III.8.4(i)]{Sc}.

(iii) If $Tf \ge f \ge0$, then $\|f\| \le \|Tf\| \le \|f\|$, so $Tf=f$ a.e. and (i) applies.
(Note that  the ergodic limit $E$ is a positive contractive projection on $F(T)$. 
We can then apply directly \cite[Theorem 4.1]{BL}, and complete the proof as  above).
\end{proof}

{\bf Remarks.} 1. For bi-stochastic p-UI Markov operators, (iii) was proved by Wu 
\cite[Corollary 3.6(b)]{Wu}.

2. Gl\"uck  \cite[Lemma 4.2]{G2} proved, without the assumption on $E$ in (i),
 that for $T$ hyperbounded as in Proposition \ref{fixed-space}, $F(T)$ is finite-dimensional .

\begin{cor} \label{ue}
Let $1 < p< \infty$ and let $T$ be a positive power-bounded operator on $L^p(S,\mu)$.
If $T$ is p-UI, then it is uniformly ergodic.
 If in addition $F(T)$ finite-dimensional,  in particular if $\|T\| \le 1$ or $T$ is hyperbounded, 
then $T$ is  quasi-compact, $\sigma(T) \cap \mathbb T$ is finite, all its points are eigenvalues   
 with finite-dimensional eigenspaces, and for some integer $d \ge 1$ all these eigenvalues 
are {\rm d}th roots of unity.
\end{cor}
\begin{proof} 
If $1 \notin \sigma(T)$, then $I-T$ is invertible and $A_n= \frac1n (I-T)^{-1}T(I-T^n) \to 0$, 
with $F(T)=\{0\}$. Otherwise, by \cite[Corollary 2.3]{G2}, the essential spectral radius
 satisfies $r_{ess}(T)<1$, so if $1 \in \sigma(T)$, then it is a pole of the resolvent, and by 
Dunford's uniform ergodic theorem \cite[p. 648]{D} $T$ is uniformly ergodic.

When $F(T)$ is finite-dimensional, which is the case if $T$ is hyperbounded by \cite{G2},
or when $\|T\| \le 1$ by Proposition \ref{fixed-space}, then
$T$ is quasi-compact by \cite{L2}. The other assertions follow from Proposition \ref{UE}.
\end{proof}


{\bf Remarks.} 1. For additional information on quasi-compact positive {\it contractions}
of $L^p$, see \cite[Theorem 2]{Sc2}.

2. Proposition \ref{no-hyper} below presents an aperiodic  Harris recurrent symmetric 
bi-stochastic Markov operator which is uniformly ergodic in $L^2$, but not hyperbounded on 
$L^2$. Corollary \ref{not-hyper} and Proposition \ref{not-hyper2} below present ergodic 
bi-stochastic Markov operators, defined by convolutions on $\mathbb T$, which (by \cite{DL}) 
are uniformly ergodic on $L^2(\mathbb T)$ (hence quasi-compact by \cite{L2}), but are not
 hyperbounded on $L^2(\mathbb T)$. Note that uniform ergodicity in $L^2$ of an
 ergodic bi-stochastic Markov operator does not imply Harris recurrence \cite{DL}.

3. If $T$ and $S$ are {\it commuting} operators on $X$, with $T$ quasi-compact and $S$ 
a contraction, then $TS$ is quasi compact.

\begin{lem} \label{product}
Let $T$ and $S$ be bounded operators on $L^p(S,\mu)$, $1 \le p< \infty$. If $T$ is
hyperbounded, so is $TS$.
\end{lem}
\medskip

\section{Limit theorems for $L^2$ quasi-compact Markov operators}

 Let $P$ be a bi-stochastic hyperbounded Markov operator  on a probability space $(S,\Sigma,\mu)$.
Recall that $P$ maps each $L^p(S,\mu)$ to itself, $1 \le p \le \infty$. 
If $P$ maps $L^1(\mu)$ to $L^q(\mu)$, then for any $1<p<q$ we have
$\|Pf\|_q \le C\|f\|_1 \le C\|f\|_p$, so $P$ maps $L^p$ to $L^q$. Similarly, if 
$P$ maps $L^p$ to $L^\infty$, then for any $p<q<\infty$ it maps $L^p$ to $L^q$.
Thus, the standing assumption for hyperbounded Markov operators is, without loss of generality,
 $1<p<q<\infty$. 
\smallskip

Let $P$ be a Markov operator with invariant probability $\mu$. The transition probability
$P(x,A)$ is not used in the definition of hyperboundedness; only the facts that $P$ 
preserves positivity, is a contraction of $L_1(\mu)$ preserving integrals, and $P1=1$ are used.
The dual operator $P^*$ satisfies the same properties (see \cite[p. 75]{F2} or \cite[p. 131]{Kr}), 
and will be called the {\it dual Markov} operator, though it need not be given by a transition 
probability (unless some  regularity assumptions on the measure space are made).

\begin{lem} \label{product-hyper}
 Let $P$ and $Q$ be   bi-stochastic Markov operators  on a probability space 
$(S,\Sigma,\mu)$. 	If $P$ is hyperbounded, so are $PQ$ and $QP$.
\end{lem}
\begin{proof}  If $P$ maps $L^p$ to $L^q$ ($q>p$), clearly also $PQ$ and $QP$ map
 $L^p$ to $L^q$.
\end{proof}

{\bf Remark.} In Proposition \ref{product=hyper} we show two commuting bi-stochastic
Markov operators which are not hyperbounded, but their product is.
\medskip

For the sake of completeness, we prove the following  fact, observed in \cite[p. 1855]{MJV}
(and  in \cite{Ha} for convolution operators).

\begin{prop} \label{all-Lp}
 Let $P$ be a hyperbounded bi-stochastic Markov operator  on a probability space 
$(S,\Sigma,\mu)$. Then $P$ and $P^*$ are hyperbounded in each $L^r$, $1<r<\infty$.
\end{prop}
\begin{proof}
Let $P$ map $L^p$ to $L^q$, $1<p<q<\infty$. Assume first that $1<r<p$. Since $P$ 
maps $L^1$ into itself and $L^p$ into $L^q$, we define $\theta \in(0,1)$ by 
$\frac1r=\theta +\frac{1-\theta}p$, and then $\frac1s=\theta  +\frac{1-\theta}q$,
and conclude from the Riesz-Thorin theorem \cite[Theorem VI.10.11]{DS} that $P$ maps 
$L^r$ to $L^s$, and $s>r$ since $q>p$.

When $p< r< \infty$, we use for the interpolation the fact that $P$ maps $L^\infty$ into
itself; then $P$ maps $L^r$ to $L^s$ where $s=rq/p$. We omit the details.

Since $P$ maps $L^p$ to $L^q$, $P^*$ is hyperbounded, mapping $L^{q'}$ to $L^{p'}$
(with $p'$ and $q'$ the dual indices). Now apply the previous part to $P^*$.
\end{proof}

{\bf Remarks.} 1. By Corollary \ref{ue} a hyperbounded Markov operator which maps $L^p$ to 
$L^q$ with $1<p<q$ is quasi-compact, hence uniformly ergodic, in $L^p$. Hence, by Proposition
\ref{all-Lp}, a hyperbounded Markov operator is quasi-compact and uniformly ergodic in all $L^r$
spaces, $1<r<\infty$. See also Corollary \ref{ros2}.

2. A hyperbounded $P$ as above need not be hyperbounded in $L^1$;
see the remarks following Proposition \ref{Lr-deriv}. 

\begin{cor}
 Convex combinations of two hyperbounded Markov operators on $(S,\Sigma,\mu)$
are hyperbounded.
\end{cor}

\begin{cor} \label{p-2}
Let $P$ be hyperbounded on $(S,\Sigma,\mu)$. Then there exists $1\le p<2$ such that $P$
maps $L^p$ to $L^2$.
\end{cor}
\begin{proof}
By Proposition \ref{all-Lp}, $P^*$ is  hyperbounded in $L^2$, so maps $L^2$ to $L^q$ for
some $q>2$. Then by duality $P$ maps $L^p$, $p=q/(q-1)$, to $L^2$.
\end{proof}

\begin{cor}
Let $P$ be bi-stochastic on $(S,\Sigma,\mu)$. The symmetrized operator 
$P_s:=\frac12(P+P^*)$ is hyperbounded if and only if $P$ is hyperbounded.
\end{cor}
\begin{proof} Let $P$ be hyperbounded. By Proposition \ref{all-Lp} both $P$ and $P^*$
are hyperbounded on $L^2$, hence so is $P_s$.

Conversely, if $P_s$ maps $L^p$ to $L^q$ with $q>p$, then $Pf\le 2P_s f$ for $f \ge 0$ yields
$$
\|Pf\|_q \le \|P|f|\,\|_q \le 2\|P_s |f|\,\|_q \le 2 \|P_s\|_{L^p \to L^q}\|f\|_p \  ,
$$
so $P$  is hyperbounded.
\end{proof}



\begin{prop} \label{ui}
Let $P$ be a bi-stochastic  Markov operator  on a probability space $(S,\Sigma,\mu)$.
If  $P$ is p-UI for some $1\le p< \infty$, then $P$ and $P^*$ are r-UI and quasi-compact 
in each $L^r$, $1<r<\infty$.
\end{prop}
\begin{proof}
For $1<r< \infty$, we apply \cite[Proposition 1.2(e)]{Wu}, with $q=1$ when $1<r<p$, 
and with $q=\infty$ when $r>p$, to conclude that $P$ is r-UI. 
 By Corollary \ref{ue}, $P$ is uniformly ergodic in $L^r$. Since by 
\cite[Corollary 3.6(b.iii)]{Wu} $F(P)$ is finite-dimensional, $P$ is quasi-compact by \cite{L2}.

Fix $1<r<\infty$. By \cite[Proposition 1.2(c) and Remark 1.3(b)]{Wu}, $P^*$ is r-UI in $L^r$,
since by the above $P$ is $r'$-UI for $r'=r/(r-1)$, and the previous part of the proof applies
 to $P^*$ on $L^r$.
\end{proof}

{\bf Remark.} Recall that a bi-stochastic $P$ may be quasi-compact in $L^2$ and not 2-UI
 \cite{BWY}, and it may be 2-UI and not hyperbounded \cite{MJV}.

\begin{prop} \label{rosenblatt}
Let $P$ be a  bi-stochastic Markov operator on a probability space $(S,\Sigma,\mu)$.
Let $Ef =\lim A_n(P)f$ for $f \in L^1(\mu)$ define the projection on the integrable invariant
functions (convergence  in $L^r$ for $f \in L^r$, $1 \le r<\infty$).
If $\|P^n-E\|_p \to 0$ for some $1\le p <\infty$, then for every $1<r<\infty$ we have
$\|P^n-E\|_r \to \infty$.
\end{prop}
\begin{proof}  Fix $1 \le r< \infty$. By the mean ergodic theorem, $Ef$ on $L^r$ projects on 
the invariant functions in $L^r$, with null space $Z_r:=\overline{(I-P)L^r}$, which is 
$P$-invariant. Let $R$ be the restriction of $P$ to $Z_r$. It is easy to see that 
$\|P^n-E\|_r \to 0$ is equivalent to $\|R^n\|_r \to 0$. We clearly have $PE=EP=E=E^2$.

Define $U=P-E$. Then $UE=EU=0$, and $U^2=PU=P(P-E)$. By induction $U^n=P^{n-1}(P-E)$,
so $\|U^n\|_r \le 2$ and also $\|U^n\|_\infty \le 2$.

For $f \in L^r$ with $Ef=0$ (i.e. $f \in Z_r$) we have $Uf=Pf =Rf$, so $U^nf =R^nf$.
Hence $\|R^n\|_r \le  \|U^n\|_r$. For any  $f\in L^r$,  $(P-E)f \in Z_r$, so
$U^n f= P^{n-1}(P-E)f = R^{n-1}(P-E)f$; hence $\|U^n\|_r\le 2\|R^{n-1}\|_r$.

Let $p<r< \infty$. We denote by $R_p$  the restriction of $P$ to $Z_p$ and by $R_r$ the
restriction to $Z_r$. Since all $L^r$ spaces are invariant under $U$, and $\|U^n\|_p \le 2$, 
$\|U^n\|_\infty \le 2$, by the Riesz-Thorin theorem \cite[Theorem (1.11), formula (1.14)]{Z}
(see also \cite[Theorem VI.10.11]{DS}), there exists $\theta \in (0,1)$ such that
$$
\|U^n\|_r \le \|U^n\|_p^\theta \|U^n\|_\infty^{1-\theta} \le
2^\theta \|R_p^{n-1}\|_p^\theta 2^{1-\theta} = \|R_p^{n-1}\|_p^\theta.
$$
 This yields $\|R_r^n\|_r \le \|U^n\|_r \le \|R_p^{n-1}\|_p^\theta  \to 0$ by the assumption, 
which is equivalent to $\|P^n-E\|_r \to 0$.

For $1<r<p$ we do the interpolation between 1 and $p$ (using $\|U^n\|_1 \le 2$).
\end{proof}

{\bf Remarks.} 1. For $P$ ergodic the theorem was proved by M. Rosenblatt 
\cite[Theorem VII.4.1]{R}; our proof is an adaptation of his. 
Note that the proof does not use positivity of $P$, nor $P1=1$, 
and applies to any contraction of $L^1(\mu)$ which is also a contraction of $L^\infty(\mu)$.

2. An example of  Rosenblatt \cite[p. 213]{R} shows that in general, even for $P$ ergodic,
 the convergence $\|P^n-E\|_r \to 0$ for every $1< r<\infty$  does not imply convergence in 
$L_1$ operator norm.

3. Rosenblatt \cite[p. 211]{R} proved also that for $P$ ergodic, $\|P^n-E\|_1 \to 0$
 is equivalent to the dual $P^*$ satisfying Doeblin's condition (so $P$ is Harris recurrent).
An example in \cite{DL} shows that $P$ ergodic with  $\|P^n-E\|_2 \to 0$ need not be
Harris recurrent.

4. Rosenblatt \cite[Lemma VII.4.1]{R} proved that for $P$ ergodic, 
the stationary Markov chain with transition probability $P(x,A)$ and initial distribution $\mu$
 is asymptotically uncorrelated if and only if $\|P^n-E\|_2 \to 0$.

\begin{cor} \label{ros2}
Let $P$ be a  bi-stochastic Markov operator on a probability space $(S,\Sigma,\mu)$.
If $P$ is uniformly ergodic in $L^p$ for some $1 \le p<\infty$, then $P$ is uniformly ergodic 
in every $L^r$, $1<r<\infty$.
\end{cor}
\begin{proof}
Put $Q=\frac12(I+P)$. By Lemma 2.1 of Foguel and Weiss \cite{FW}, $\|Q^n(I-Q)\|_p \to 0$; 
since $P$ and $Q$ have the same invariants, uniform ergodicity in $L^P$ yields 
$\|Q^n -E\|_p \to 0$. By Proposition \ref{rosenblatt}, $\|Q^n-E\|_r \to 0$ for any
$1<r< \infty$. Then $(I-Q)L^r=(I-P)L^r$ is closed in $L^r$, so $P$ is uniformly ergodic in $L^r$.
\end{proof}

{\bf Remarks.} 1. When $P$ is also ergodic, the corollary was proved in \cite[Corollary 3.5]{DL}.

2. By Corollary \ref{ros2}, a bi-stochastic  $P$ is uniformly ergodic in some $L^p(\mu)$,
 $1<p<\infty$, if and only if it is {\it $L^2$-uniformly ergodic}.
\smallskip

{\bf Definition.} For a bi-stochastic Markov operator quasi-compact in $L^p(\mu)$, the 
value $d$ in Proposition \ref{qc} is called the {\it period} of $P$ (in $L^p$); we call $P$ 
{\it periodic} if $d>1$, and {\it aperiodic} when $d=1$.

\begin{prop} \label{qc-Lp}
Let $P$ be a bi-stochastic Markov operator on a probability space $(S,\Sigma,\mu)$.
Let $Ef =\lim A_n(P)f$ for $f \in L^1(\mu)$ define the projection on the integrable invariant
functions.
If $P$ is quasi-compact in $L^p$ for some $1<p< \infty$, then it is quasi-compact in every
$L^r$, $1 < r<\infty$, and the period of $P$ in $L^r$ equals its period in $L^p$
(we call this common period the {\rm period} of $P$).
\end{prop}
\begin{proof} Since $\mu$ is an invariant probability, $P$ is conservative (e.g. 
\cite[p. 117]{Kr}), and all integrable invariant functions are measurable with respect to the 
$\sigma$-algebra of invariant sets $\Sigma_i:= \{A\in \Sigma: P1_A=1_A\}$, since the limit 
$Ef$ in the ergodic theorem is the conditional expectation $E(f|\Sigma_i)$
(e.g. \cite[p. 128]{Kr}, or \cite[p. 80]{F2}). By quasi-compactness in $L^p$, 
$F(P) \cap L^p$ is finite-dimensional, so $L^p(S,\Sigma_i,\mu)$ is finite-dimensional. 
This implies that $\Sigma_i$ is finite modulo $\mu$, so  for $1<r< \infty$ also 
$F(P)\cap L^r$ is finite-dimensional. By Corollary \ref{ros2} $P$ is uniformly ergodic 
in $L^r$, so by \cite{L2}, $P$ is quasi-compact on $L^r$.
\smallskip

By the above, $P$ is quasi-compact in $L^2$, and let $d$ be its period in $L^2$. 
Fix $1<r< 2$, and let $d_r$ be the period of $P$ in $L^r$. When $f \in L^2$ satisfies
$Pf = \lambda f$ with $ \lambda \in \mathbb T$, then $\lambda^d=1$, and since 
$f \in L^r$, also $\lambda^{d_r}=1$. Minimality of the period yields $d \le d_r$.

The dual of a quasi-compact operator is quasi-compact in the dual space, hence the dual
 Markov operator $P^*$ is quasi-compact in $L^r$, $1<r<\infty$, with the period of $P^*$
in $L^r$ the same as the period of $P$ in the dual $L^{r'}$ ($r'=r/(r-1)$).
Now fix $2<r<\infty$, with $d_r$ the period of $P$ in $L^r$. The above argument
yields $d_r \le d$.But $r'<2$, so $d \le d_{r'}=d_{r}$, hence $d=d_r$ for $r>2$, 
and by duality $d=d_r$ also for $1<r<2$.
\end{proof}

\begin{cor} \label{invariants}
Let $P$ be a bi-stochastic Markov operator on $(S,\Sigma,\mu)$. If $P$ is quasi-compact 
in $L^2$, then every integrable eigenfunction corresponding to a unimodular eigenvalue,
in particular every integrable invariant function, is bounded.
\end{cor}
\begin{proof}
We saw that $\Sigma_i$ is finite, so generated by finitely many atoms. If $Pf=f \in L^1$,
then $f$ is $\Sigma_i$-measurable, so bounded.

 If $Pf=\lambda f$ with $|\lambda|=1$ and $f \in L^1(\mu)$, then $P|f| \ge |Pf|= |f|$, and 
by invariance of $\mu$ we have $P|f|=|f|$;  hence $f$ is bounded.
\end{proof}

{\bf Remark.}  In general, the dimension of the eigenspace of a unimodular eigenvalue of
a bi-stochastic Markov operator $P$, which is mean ergodic, is not more than the dimension of 
$F(P)$, by \cite[Theorem 2]{L2}. In particular, if $P$ is ergodic, then the eigenspaces of unimodular 
eigenvalues are one-dimensional, and the eigenfunctions have constant absolute value.

\begin{cor} \label{all-r}
Let $P$ be a hyperbounded bi-stochastic Markov operator on $(S,\Sigma,\mu)$.
Then for every $1<r<\infty$, $P$ is quasi-compact in $L^r$.
\end{cor}
\begin{proof} Use Corollary \ref{ue} and Proposition  \ref{qc-Lp} when $P$ is 
hyperbounded on $L^p$.
\end{proof}
%
%
%
%

\begin{prop} \label{exponential}
Let $P$ be a bi-stochastic Markov operator on a probability space $(S,\Sigma,\mu)$,
and assume that for some $1<p<\infty$,  $P$ is quasi-compact in $L^p(\mu)$ with period $d$.
 Let $E_d := \lim A_n(P^d)$ be the projection on $F(P^d)$. Then:

(i) for any $1<r< \infty$ there exist $C_r>0$ and $\rho_r <1$ such that 
  $\|P^{nd}-E_d\|_r \le C_r\rho_r^n$.

(ii) If  $f \in L^r(\mu)$, $1<r< \infty$, then 
$\lim_{n\to\infty} P^{nd} f=E_d f$ a.e. and in $L^r$.
In particular, if $P$ is aperiodic, $P^nf \to Ef$ a.e.

(iii) If $f \in L^1(\mu)$, then  $\|P^{nd}f -E_df\|_1 \to 0$.
\end{prop}
\begin{proof}
By Proposition \ref{qc-Lp} $P$ is quasi-compact in any $L^r$, $1<r<\infty$, with the same
period $d$. By Proposition \ref{qc}, $P^{nd}$ converges in $L^r$ operator norm,
necessarily to $E_d$ by the mean ergodic theorem.
 We now apply Proposition 3.1 of \cite{DL} to obtain the exponential rate (i).


Since $P^d E_d=E_d$, for $1<r< \infty$ and  $f \in L^r(\mu)$  (i) yields 
$$\sum_{n=1}^\infty \|P^{nd}f -E_df\|_1 \le \sum_{n=1}^\infty \|P^{nd}f -E_d f\|_r < \infty,$$
 so by Beppo Levi $\sum_{n=1}^\infty |P^{nd}f -E_d f| < \infty$ a.e.,  which implies (ii).

Since $\|P^{nd}f -E_df\|_1 \to 0$ for $f \in L^2$, say, (iii) follows by continuity.
\end{proof}

\begin{cor} \label{doeblin}
Let $P$ be an ergodic hyperbounded bi-stochastic Markov operator with period $d$.
If $P$ maps $L^1(\mu)$ to $L^q(\mu)$ for some $q>1$, then:

(i) There exist $C_1>0$ and $\rho_1 <1$ such that $\|P^{nd}-E_d\|_1 \le C_1 \rho_1^n$.

(ii) $P$ is uniformly ergodic in $L^1$.

(iii) For any $f \in L^1(\mu)$, $\ \lim_{n\to\infty} P^{nd}f =E_df$ a.e. and in $L^1$. In particular,
if $d=1$,  then $\lim_{n\to\infty} P^{n} f=\int f\,d\mu$ a.e. for every $f \in L^1$.
\end{cor}
\begin{proof}
We can assume $q< \infty$. Fix some $p \in(1,q)$. Then $P$ maps also $L^p(\mu)$ into 
$L^q(\mu)$; also $P^d$ maps $L^1$ and $L^p$ into $L^q$, and put 
$C:=\|P^d\|_{L^1\to L^q}$. By Corollary \ref{ue}, $P$ is quasi-compact in $L^p(\mu)$.

(i) Fix $f \in L^1$. We use $E_dP^d=P^dE_d=E_d=E_d^2$ and obtain
$$
\|P^{nd}f-E_df\|_1 \le \|P^{nd}f -E_df\|_q = \|P^{(n-1)d} -E_d\|_q\|P^d(f-E_df)\|_q \le
$$
$$
C \|P^{(n-1)d} -E_d\|_q\|(f-E_df)\|_1 \le 2 C \|P^{(n-1)d} -E_d\|_q\|f\|_1.
$$
Hence $\|P^{nd}-E_d\|_1 \le 2C \|P^{(n-1)d}-E_d\|_q \le 2C\cdot C_q\rho_q^{n-1}$, 
by Proposition \ref{exponential}(i).

(ii) follows from (i).

(iii) follows from (i) and Beppo Levi's theorem.
\end{proof}

{\bf Remark.} Theorem \ref{gluck} does not apply to hyperbounded operators in $L^1$,
 so does not yield directly (ii) of the corollary.
\medskip

{\bf Example 1.} {\it Hyperbouned Markov operators}

 On $(S,\mu)$ we define a Markov operator by $Pf(x)= \int k(x,y)f(y)d\mu(y)$,
with a kernel $k(x,y) \ge 0$ satisfying $\int k(x,y)d\mu(y)=1$ for a.e. $x$ and 
$\int k(x,y)d\mu(x)=1$ for a.e. $y$. Then $P$ is bi-stochastic, 
hence a contraction of every $L^p(\mu)$. $p\ge 1$. If $k(x,y)$ is bounded, then 
$PL^1(\mu) \subset L^\infty(\mu)$, so $PL^p(\mu) \subset L^q(\mu)$ for any $1<p<q< \infty$;
by \cite{Wu} $P$ is 1-UI (i.e. weakly compact in $L^1$; see also \cite[Exercise VI.9.57]{DS}).
In fact, if for some $q>1$ and $M<\infty$ we have $\int |k(x,y)|^q d\mu(y) \le M$ for
a.e. $x$, then $P$ maps $L^1$ to $L^q$, with $\|P\|_{L^1\to L^q} \le M$ 
\cite[Exercise VI.9.59]{DS}; hence $P$ is 1-UI.

If  $\int \int |k(x,y)|^q d\mu(x)d\mu(y) < \infty$ for some $q>2$,
 then $PL^2(\mu) \subset L^q(\mu)$ \cite[p. 480]{CM}; in fact, for 
$q':= \frac q{q-1} < 2$ we have $PL^{q'} \subset L^q$. By Proposition \ref{all-Lp},
$P$ is hyperbounded in every $L^r$, $1<r<\infty$.
\medskip

{\bf Remarks.} 1.  The Markov operator $P$  defined in Example 1 is Harris recurrent (by 
the analytic definition in \cite[p. 492]{F3}; for the probabilistic definition see \cite{MT}).

2. When $P$ is an ergodic hyperbounded Markov operator mapping $L^1(\mu)$
to $L^q(\mu)$, then $P$ is given by a kernel (see \cite[Exercise VI.9.59]{DS}); hence
$P$ is Harris recurrent. If in addition $P$ is aperiodic, then $P^*$ satisfies
Doeblin's condition, by Corollary \ref{doeblin}(ii) and  \cite[p. 212]{R}.

3. Let $P(x,A)$ define  an ergodic hyperbounded Markov operator $P$ on $(S,\Sigma,\mu)$, 
and let $\{\xi_n\}$ be the stationary Markov chain with transition probability $P(x,A)$
and initial distribution $\mu$. By Corollary \ref{all-r}, $P$ is uniformly ergodic in $L^2(\mu)$,  
so for every $f \in L^2(S,\mu)$ with  $\int f\,d\mu=0$ the sequence 
$\{ \frac1{\sqrt n} \sum_{k=1}^n f(\xi_k)\}$ satisfies the central limit theorem (CLT), 
by a result of Gordin and Lifshitz, see \cite{DL}. 
Davydov \cite{Da} constructed an aperiodic positive recurrent Markov chain with state space 
$\mathbb Z$, for which the above the CLT fails for some $f \in L^2(\mu)$ with integral zero. 
The operator given by Davydov's transition matrix is Harris and not uniformly ergodic
in $L^2$, hence not hyperbounded. $P^n$ converges in the $L^2$-strong operator
 topology, but not in $L^2$-operator norm.

4. Let $\{\xi_n\}$ be the stationary Markov chain defined by an {\it aperiodic}
ergodic 2-UI (in particular hyperbounded)  $P(x,A)$ on $(S,\Sigma,\mu)$. 
By Proposition \ref{exponential},  $\|P^n-E\|_2 \to 0$. In this case, 
El Machkouri et al. \cite[Theorem 3.1]{MJV} proved a limit theorem for the 
distribution of $\{\frac1{B_n}\sum_{k=1}^n f(\xi_k)\}$ 
for appropriate $B_n$, when $f$ on $S$ has no variance but has heavy tails.

5. A hyperbounded bi-stochastic Markov operator need not be Harris recurrent. 
An example is given in Theorem \ref{riesz}.
\medskip

 Under the assumptions of Theorem \ref{gluck}, if  $\{T^n\}$ converges weakly, 
then (by Corollary \ref{ue}) it converges in operator norm ($d=1$); however, in general 
this is not so. 
\medskip

{\bf Example 2.} {\it A periodic symmetric ergodic  hyperbounded Markov operator}

Let $S=[0,1]$ with $\mu$ the Lebesgue measure, and define a bounded kernel by $k(x,y)=2$ 
on  $([0,\frac12] \times [\frac12,1]) \cup ([\frac12,1] \times [0,\frac12])$ and zero elsewhere 
on $[0,1]\times [0,1]$. Then the corresponding  (ergodic and symmetric) Markov operator $P$ 
as defined in example 1 has $Pf = -f$ for $f=1_{[0.\frac12]} -1_{[\frac12,1]}$. Note that
"spectral gap" in the sense of Miclo \cite{M} is only "spectral gap {\it near} 1" (as in \cite{G2}).

\begin{prop} \label{no-hyper}
There exists an aperiodic Harris recurrent symmetric Markov operator $P$ with $P^n$ 
convergent in $L^2$-operator norm, which is not hyperbounded.
\end{prop}
\begin{proof} Denote by $Q$ the operator of Example 2, and define $P=\frac12(I+Q)$.
Clearly $Q$ is  ergodic, and since $I-P=\frac12(I-Q)$, also $P$ is ergodic; symmetry of $Q$
implies that of $P$. Since $Q$ is defined by a kernel, $P$ is Harris recurrent. 
By the construction of Example 1 $Q$ maps $L^p$ into $L^q$ for any $1 \le p<q\le \infty$, 
so it is uniformly ergodic in every $L^p(\mu)$, $1<p<\infty$, by Corollary \ref{ue}. 
By Foguel and Weiss \cite[Lemma 2.1]{FW}, $\|P^n(I-P)\|_p \to 0$, so with the uniform 
ergodicity $P^n$ converges in $L^p$ operator norm ($1<p<\infty$).

If $P$ were hyperbounded, it would map $L^p(\mu)$ to $L^q(\mu)$, for some $1<p<q< \infty$.
Let $f \in L^p$ which is not in $L^q$. Since $Q$ maps $L^p$ into $L^q$ (see Example 1), 
$\frac12(I+Q)f=Pf \in L^q$ implies $f \in L^q$, a contradiction. Hence $P$ is not hyperbounded.
\end{proof}
\smallskip

\medskip

\section{Cyclic behavior of ergodic bi-stochastic Markov operators}

In this section $P$ is an {\it ergodic} Markov operator on $(S,\Sigma,\mu)$ with $\mu$
an invariant {\it probability}. Then $P$ is conservative, and by Hopf's pointwise ergodic theorem
\cite[Theorem VIII.6.6]{DS},\cite[p. 80]{F2}, \cite[Theorem 1.7.2]{Kr} we have a.e. convergence
of the averages $A_n(P)f(x)$ for any $f \in L^1(\mu)$.

 On the other hand, let $P$ be a transition probability and  $\tilde \mu$ a probability 
on $(S,\Sigma)$ such that $Pf=0$ $ \tilde \mu$ a.e whenever $f=0\ $ $\tilde\mu$ a.e. 
Then $L^\infty(\tilde\mu)$ is invariant under $P$, and $P$ on $L^\infty(\tilde\mu)$  is the 
dual of a positive contraction on $L^1(\tilde\mu)$. If $P$ on $L^\infty(\tilde \mu)$ is
 conservative, then a.e. convergence of the averages $A_n(P)f(x)$ 
(or even convergence of their integrals with respect to $\tilde\mu$) for every $f \in L^\infty$, 
implies the existence of a probability $\mu \sim \tilde\mu$ which is invariant \cite{LS}.
 This justifies our assumption that $P$ is bi-stochastic on $(S,\Sigma,\mu)$ for studying the
convergence of the iterates $P^n$.
 \smallskip

When $P$ is ergodic and bi-stochastic, it is quasi-compact in $L^p(\mu)$ if and only if it is
uniformly ergodic in $L^p(\mu)$, by \cite{L2}. When $P$ is quasi-compact in $L^1(\mu)$ and 
aperiodic, then $P^*$ satisfies Doeblin's condition \cite[p. 211]{R}, hence both $P^*$ and $P$
are Harris recurrent. However, $L^2$-quasi-compactness does not imply Harris recurrence
\cite{DL}. Even hyperboundedness of $P$ does not imply Harris recurrence  (Theorem
\ref{riesz} below).
\smallskip

{\bf Remark.} For $P$ bi-stochastic on $(S,\Sigma,\mu)$  (without assuming
ergodicity), quasi-compactness in $L^p(\mu)$ for some $1\le p < \infty$ implies that the
$\sigma$-algebra of invariant sets $\Sigma_i$ is finite (see proof of Proposition \ref{qc-Lp});
hence the reduction to each of the finitely many atoms of $\Sigma_i$ will be an {\it ergodic}
bi-stochastic Markov operator as above. 
\smallskip

{\bf Definition.} A power-bounded operator $T$ on a Banach space $L$ is called 
{\it constrictive} if there exists a compact set $\mathcal K \subset L$ such that
\begin{equation} \label{constrict}
\text{dist}(T^n x,\mathcal K) \underset{n\to \infty}\to  0 \quad \text{ for every } \ \|x\| \le 1.
\end{equation}
By the definition, for every $x \in L$ the orbit $(T^nx)_{n\ge 1}$ is precompact, so $T$ is
strongly almost periodic, hence mean ergodic.

\smallskip

$L^1$-constrictive bi-stochastic Markov operators on a probability space $(S,\Sigma,\mu)$
were introduced and studied by Lasota, Li and Yorke \cite{LLY}; for $L=L^1(\mu)$  they proved   
\begin{equation} \label{lly}
\lim_{n\to\infty} \|T^n(x-\sum_{j=1}^r \varphi_j(x)y_j)\| \to 0 \quad \text{for every } \ x \in L,
\end{equation}
with $y_j$ non-negative unit vectors with disjoint supports, and $T$ permutes the 
$(y_j)_{1\le j \le r}$.  It follows that for some $d \le r!$, $\ T^{nd}x$ converges strongly 
to $\sum_{j=1}^r \varphi_j(x)y_j$.
Obviously, if $P$ is ergodic with $P^nf \to \int f\,d\mu$ in $L^1(\mu)$, then it is constrictive. 
Komorn\'\i k \cite{Kom} proved \eqref{lly} when the
 Markov operator is only {\it weakly} constrictive in $L^1$, i.e.  \eqref{constrict} holds 
with $\mathcal K$ only weakly compact.

Bartoszek \cite[Theorem 2]{Ba} proved that {\it if $T$ is a quasi-compact positive 
contraction on a Banach lattice $L$, then $T$ is constrictive}, with the convergence in 
\eqref{constrict} {\it uniform over the unit ball}. In \cite[Theorem 1]{Ba} he proved that  
if $T$ is a constrictive positive contraction on $L$, then there exist $r$ positive unit 
vectors $y_1,\dots,y_r$ and $r$ positive functionals $\varphi_1\dots,\varphi_d$ in $ L^*$ 
such that \eqref{lly} holds, and $T$ permutes the $(y_j)_{1\le j\le r}$. 

Sine \cite{Si} used the deLeeuw-Glicksberg decomposition (see, e.g. \cite{Kr}) to study 
general constrictive contractions in Banach spaces.
\medskip

The next results are inspired by the work of Foguel \cite{F1},\cite{F2},\cite{F5},\cite{F3},\cite{F6}.
\medskip


{\bf Definition.} The {\it deterministic $\sigma$-algebra} of a Markov operator $P$ is defined by
$$
\Sigma_D=\Sigma_D(P):= 
\{A\in \Sigma: \text{ for every } n\ge 1, \quad P^n 1_A=1_{A_n} \text{ for some }\ A_n \in \Sigma\}.
$$

\noindent
The proof that $\Sigma_D$ is a $\sigma$-algebra is in \cite[p. 7]{F2}. In general,
$\Sigma_D(P^*) \ne \Sigma_D(P)$ \cite[p. 78]{F2}; in that example, $P^*$ is constrictive
and $P$ is not. If $P$ is Harris recurrent, then $\Sigma_D$ is atomic \cite[p. 58]{F2} 
(proof corrected in \cite{F4}). 

It is shown in \cite[p. 106]{R} (see also \cite[p. 87]{F2}) that if $P$ is  bi-stochastic on 
$(S,\Sigma,\mu)$, then 
$$
L^2(S,\Sigma_D(P),\mu) =\{g \in L^2(\mu): \|P^n g\|_2=\|g\|_2 \text{ for every } n\in\mathbb N\}.
$$

\noindent
It follows from \cite[Theorem A, p. 85]{F2} that if $f \in L^2$ satisfies $E (f|\Sigma_D)=0$, 
then $P^nf \to 0$ weakly in $L^2$. Consequently, if 
$E(f|\Sigma_D)=0$ and for some subsequence $P^{n_k}f$ converges in $L^2$-norm, then  
$\|P^nf\|_2 \to 0$. In particular, If $P$ is $L^1$-constrictive, then $\|P^nf\|_2 \to 0$
whenever  $E(f|\Sigma_D)=0$. In general, $\|P^n f\|_2 \to 0$ implies $E(f|\Sigma_D)=0$
\cite[p. 108]{R}. However, {\it $\|P^n f\|_2 \to 0$ for every $f$ with $E(f|\Sigma_D)=0$ 
if and only if the strong limit $\lim_{k \to \infty} P^{*k}P^k$ (which always exists -- 
\cite[p. 108]{R}) is a projection}  \cite[Lemma 3, p. 108]{R}, which is $E(f|\Sigma_D)$; an example 
by Rosenblatt  \cite[p. 113]{R}, with $\Sigma_D$ trivial (hence all  powers ergodic, by Lemma 
\ref{k-invariants} below), shows that the above strong convergence need not hold in general. 

\begin{lem} \label{k-invariants}
Let $P$ be an ergodic bi-stochastic Markov operator on $(S,\Sigma,\mu)$.
Then for every $k \ge 1$ we have 
$$
\Sigma_{i,k}:= \{A \in \Sigma:\, P^k 1_A=1_A\} \subset \Sigma_D(P).
$$
\end{lem}

{\bf Remarks.} 1. The lemma was first proved in \cite[Lemma 2.1.8]{F5}. An accessible proof is
in \cite[Lemma 1.2]{F6}.

2. The example of $P$ induced by an irrational rotation of the unit circle 
$\mathbb T$ shows that in general we may have all powers of $P$ ergodic, i.e. 
$\Sigma_{i,k}$ trivial  for every $k$, while $\Sigma_D=\Sigma$.

\begin{prop} \label{foguel}
Let $P$ be an ergodic bi-stochastic Markov operator on $(S,\Sigma,\mu)$. Then $\Sigma_{i,k}$
is finite for any $k > 1$, and has at most $k$ atoms.

Moreover, for fixed $k>1$ there exist $d|k$ and atoms $A_0,\dots,A_{d-1}$ of $\Sigma_{i,k}$
which are disjoint, generate $\Sigma_{i,k}$, and $P1_{A_j}=1_{A_{j+1}}$
for $0\le j<d$ (with $A_d=A_0$).
\end{prop}
\begin{proof} Fix $k$. If $\Sigma_{i,k}$ is not trivial, let $A \in \Sigma_{i,k}$ with 
$0<\mu(A)<1$.  By Lemma \ref{k-invariants}
$A \in \Sigma_D$, so there are $B_j,\  j=0,1,\dots,k-1$ with $P^j 1_A=1_{B_j}$. Then 
$\mu(B_j)=\mu(A)$ and by ergodicity $1_A \le \sum_{j=0}^{k-1}P^j 1_A = k\mu(A)$.
Hence $\mu(A) \ge 1/k$ for any $A \in \Sigma_{i,k}$ with $\mu(A)>0$. This implies that
$\Sigma_{i,k}$ is atomic, with at most $k$ different (hence disjoint) atoms.
\smallskip

Now fix $k>1$, and let $A$ be an atom of $\Sigma_{i,k}$. By Lemma \ref{k-invariants},
$P1_A=1_{A_1}$, and clearly $A_1 \in \Sigma_{i,k}$, with $\mu(A_1)=\mu(A)>0$. 
We show that $A_1$ is an atom of $\Sigma_{i,k}$. Let $B \subset A_1$ with $\mu(B)>0$
be in $\Sigma_{i.k}$. Then  $P^{k-1}1_B=1_C$ for some $C \in \Sigma_{i,k}$, and
$P^{k-1}1_B \subset P^{k-1}1_{A_1}=1_A$. But then $C \subset A$ with $\mu(C)=\mu(B)>0$,
so $C=A$, and $1_B=P1_C=1_{A_1}$. Hence $A_1$ is an atom of $\Sigma_{i,k}$. 

Let $d$ be smallest integer with $P^d1_A=1_A$, and set $1_{A_j}=P^j1_A$ for $0\le j<d$.
By the above, the $A_j$ are atoms, and disjoint by minimality of $d$. Since $P^k1_A=1_A$,
the minimality of $d$ yields $d|k$. By definition, $P1_{A_j}=1_{A_{j+1}}$. By disjointness
and ergodicity, $1_A \le \sum_{j=0}^{k-1}P^j1_A =1$, so $\cup_{0\le j<d}A_j =S$.
Finally, if $B \in \Sigma_{i,k}$ then $ B\cap A_j$ is $A_j$ or a null set, so $A_0\dots,A_{d-1}$
generate $\Sigma_{i,k}$.
\end{proof}

The disjoint sets $A_0,\dots,A_{d-1}$ obtained in Proposition \ref{foguel} are 
{\it cyclically moved} by $P$, with period $d$, and form a {\it cycle}. 
Non-trivial cycles exist when $P^k$ is not ergodic for some $k>1$.
If $f \in L^1(\mu)$ is supported in $A_0$, then $P^jf $ is supported in $A_j$.
Note that Harris recurrence is not assumed.

Note that if $A_0,\dots,A_{k-1}$ are disjoint with
$P1_{A_j}=1_{A_{j+1}}$ for $0\le j \le k-1$ (with $A_k=A_0$), then $P^k1_{A_j}=1_{A_j}$,
so $A_j \in \Sigma_{i,k}$ for $0\le j \le k-1$.
\smallskip

{\bf Remark.} Foguel \cite[ Theorem 2.1.10]{F5} proved Proposition \ref{foguel} 
for $P$ conservative and ergodic, without the assumption of  an invariant probability.
Since \cite{F5} is not readily available, we have included the proof for our situation.

\begin{cor} \label{Lly}
Let $P$ be an ergodic bi-stochastic Markov operator on $(S,\Sigma,\mu)$. Then $P$ is
$L^1$-constrictive if and only if for some $d \ge 1\ $ $P^{nd}$ converges in the strong 
operator topology of $L^2(\mu)$, as $n \to \infty$.
\end{cor}
\begin{proof} If $P$ is $L^1$-constrictive (or even $L^1$-weakly constrictive), 
$L^1$-strong convergence of $P^{nd}f$ follows from \eqref{lly}. 
Since $P$ contracts $L^\infty$-norms, for bounded $f$ we have $L^2$ convergence
of $P^{nd}f$, hence $P^{nd}$ converges in the $L^2$ strong operator topology.

Conversely, if $P^{nd}$ converges strongly in $L^2$, it also does in $L^1$. The limit is a 
projection on the integrable $P^d$-invariant functions, which are $\Sigma_{i,d}$-measurable.
By Proposition \ref{foguel}, $L^1(\Sigma_{i,d},\mu)$ is finite-dimensional, so its unit ball 
is compact in $L^1(\mu)$; hence $P$ is $L^1$-constrictive.
\end{proof}

{\bf Remark.} The result \eqref{lly}, proved for $L^1$ (weakly) constrictive bi-stochastic
Markov operators in \cite{LLY} and in \cite{Kom}, is used only in Corollary \ref{Lly}.

\begin{prop} \label{deterministic} 
Let $P$ be an ergodic bi-stochastic Markov operator on $(S,\Sigma,\mu)$.  If  $k \ge 1\ $ 
is an integer such that $P^{nk}$ converges strongly in $L^2$ as $n \to \infty$, then
$$
\Sigma_D(P) = \Sigma_{i,k}:= \{A \in \Sigma:\, P^k 1_A=1_A\}.
$$
\end{prop}
\begin{proof}
Every $P^k$-invariant integrable function is $\Sigma_{i,k}$-measurable, and by the ergodic
theorem $E_k f := \lim_n A_n(P^k)f = E(f|\Sigma_{i,k})$ for $f \in L^1$. By assumption,
$P^{nk}f \to E_kf$ also  in $L^1$.

Let $A \in \Sigma_D$, with $P^n 1_A=1_{A_n}$. Then
$1_{A_{nk}} = P^{nk}1_A \to E_k 1_A$  in $L^2$-norm.
Now
$$
\mu(A)= \int P^n1_A d\mu= \int 1_{A_n}d\mu= \|1_{A_n}\|_2^2 =\|P^n1_A\|_2^2 \le
\|1_A\|_2^2=\mu(A).
$$
Hence $\|E_k 1_A\|_2^2  =\lim_n \|P^{nk}1_A\|_2^2= \mu(A)= \|1_A\|_2^2.$
Since $E_k$ is an orthogonal projection, this equality implies $E_k 1_A=1_A$,
 i.e. $A \in \Sigma_{i,k}$. Thus $\Sigma_D \subset \Sigma_{i,k}$. 

The reverse inclusion holds by Lemma \ref{k-invariants}, so $\Sigma_D = \Sigma_{i,k}$.
\end{proof}

{\bf Remarks.} 1. When $k=1$, $\Sigma_{i,1}=\Sigma_i = \Sigma_D$, so in the "complete 
mixing" case $\Sigma_D(P) $ is trivial. In the general ergodic case, strong convergence  in 
$L^2$ of $P^{nk}$ implies that $\Sigma_D(P)$ is atomic with at most $k$ atoms. 
For additional information see \cite{F5} and \cite{KL}.

2. The example in \cite[p. 113]{R} shows that $\Sigma_{i,k}=\Sigma_D$ for every $k$
 does not imply strong convergence of $(P^{nd})_{n \ge 1}$ for any $d$.

3. The example in \cite[p. 78]{F2} shows that $P$ in the proposition need not be Harris.

\begin{prop} \label{divisor}
Let $P$ be an ergodic bi-stochastic Markov operator on $(S,\Sigma,\mu)$. Assume that
for some integer $k_0 \ge 1\ $ $P^{nk_0}$ converges strongly in $L^2$ as $n \to \infty$,
and let $d_0$ be the smallest such $k_0$.  Then 
	$(P^{nk})_{n \ge 1}$ converges strongly in $L^2$ if and only if $d_0|k$.
\end{prop}
\begin{proof} 
Obviously $(P^{nk})_{n\ge 1}$ converges strongly for $k=md_0$.

 We have to prove the converse only when $d_0 >1$. 
Let $d \le d_0$ be the number of disjoint atoms in $\Sigma_{i,d_0}$, given by Proposition 
\ref{foguel}, which  are cyclically moved by $P$, with $P^d1_{A_j}=1_{A_j}$ for the atoms 
$A_0,\dots,A_{d-1}$. Since $d|d_0$ and the atoms generate $\Sigma_{i,d_0}$, 
we have $\Sigma_{i,d}=\Sigma_{i,d_0}$. 
With our previous notations, $E_d=E_{d_0}$, and $P^d E_{d_0}=P^d E_d=E_d$.

We now prove $d=d_0$. Since $d|d_0$, we write $d_0=md$. Fix $f \in L^2$ and $\varepsilon > 0$. 
For $n\ge N$ we have $\|P^{nd_0}f -E_{d_0}f\| < \varepsilon$. Hence for $k \ge 1$ we have
$$
\|P^{(Nm+k)d}f -E_df\| \le \|P^{kd}\|\cdot \|P^{Nmd}f -E_{d_0}f\| \le
 \|P^{Nd_0}f -E_{d_0}f\| < \varepsilon.
$$
This proves that $(P^{nd})_{n\ge 1}$ converges strongly, and minimality of $d_0$ yields $d=d_0$. 

Assume $(P^{nk})_{n\ge 1}$ converges strongly. By Proposition \ref{deterministic}, 
$\Sigma_{i,k}=\Sigma_D=\Sigma_{i,d_0}$, and $d= d_0 \le k$. Let $k=\ell \mod d$ with 
$0\le \ell <d$.  For the atom $A_0$, which is in $\Sigma_D$, we have 
$P^\ell 1_{A_0}=P^k1_{A_0}=1_{A_0}$. 
By construction, $P^j1_{A_0} \ne 1_{A_0}$ for $1 \le j \le d-1$, so $\ell=0$, which means $d|k$.
\end{proof}

\begin{cor}
Let $P$ be an ergodic bi-stochastic Markov operator on $(S,\Sigma,\mu)$. Assume that
 for some integer $k_0 \ge 1\ $ $P^{nk_0}$ converges strongly in $L^2$ as $n \to \infty$,
and let $ k \in \mathbb N$. Then $(P^{nk})_{n\ge 1}$ converges strongly in $L^2$ if
and only if $\Sigma_{i,k}=\Sigma_D(P)$.
\end{cor}
\begin{proof} Assume that $\Sigma_{i,k}=\Sigma_D(P)$, and let $d=d_0$ as defined in
Proposition \ref{divisor}. The last three sentences of the proof of Proposition \ref{divisor} show 
that $d|k$; hence Proposition  \ref{divisor} yields that  $(P^{nk})_{n\ge 1}$ 
converges strongly in $L^2$.

The converse is in Proposition \ref{deterministic}.
\end{proof}

{\bf Remark.} By Corollary \ref{Lly}, the assumption in Proposition \ref{divisor} and its 
corollary is that $P$ is ergodic, bi-stochastic and $L^1$-constrictive.

\begin{theo} \label{periodic-convergence}
Let $P$ be an ergodic bi-stochastic Markov operator on $(S,\Sigma,\mu)$. Assume that
 for some integer $k \ge 1\ $ $P^{nk}$ converges strongly in $L^2$ as $n \to \infty$,
and let $d>1$ be the smallest such $k$. Then:

(i)  $\Sigma_{i,d}=\Sigma_D(P)$, $\Sigma_D$ is atomic, generated by $d$ disjoint  
atoms $A_0,\dots,A_{d-1}$, which satisfy $P 1_{A_j}=1_{A_{j+1}}$.

(ii) The unimodular eigenvalues of $P$ in $L^1(\mu)$ are precisely
all {\rm d}th roots of unity, the corresponding eigenspaces are one-dimensional,
and the corresponding eigenfunctions have constant absolute value.

(iii)  For every $f \in L^p(\mu)$, $1 \le p< \infty$, and $0 \le j \le d-1$ we have
\begin{equation} \label{L1-limits}
\lim_{n\to\infty}\big\|P^{nd+j}f - 
d\sum_{\ell=0}^{d-1} \big(\int_{A_\ell}f \,d\mu\big) 1_{A_\ell\dot{+}j} \big\|_p = 0,
\end{equation}
where $\ell\dot{+}j$ is addition modulo $d$.

(iv)  For every $f \in L^p(\mu)$, $1 \le p< \infty$ we have
\begin{equation} \label{conditional}
\lim_{n\to\infty}\big\|P^n \Big( f -
d \sum_{\ell=0}^{d-1} \big( \int_{A_\ell}f \,d\mu\big) 1_{A_\ell}\Big) \big\|_p = 0.
\end{equation}
\end{theo}
\begin{proof}
(i) The equality $\Sigma_{i,d}=\Sigma_D$ and the existence and properties of the atoms of
$\Sigma_{i,d}$ follow form the previous results. Invariance of $\mu$ yields $\mu(A_j)=\mu(A_0)$, 
so $\mu(A_j)=d^{-1}$.

(ii) If $\lambda$ is a unimodular eigenvalue, then by assumption $\lambda^{nd}$ converges as 
$n \to \infty$, hence $\lambda^d=1$.
 Conversely, for $\lambda^d=1$  define $f= \sum_{j=0}^{d-1} \bar\lambda^j 1_{A_j}$;
then $Pf=\lambda f$, since $P1_{A_{d-1}}=1_{A_0}$.
The last statement in (ii) follows  by ergodicity.

(iii) Note that for strong  convergence of $P^n$ the limit is $\int f\ d\mu$, by ergodicity, 
so we prove \eqref{L1-limits}  for $d >1$. The convergence \eqref{L1-limits} in $L^1$ for all 
integrable functions easily implies $L^p$-norm convergence for $L^p$ functions, so we prove 
\eqref{L1-limits} for $p=1$.

For $j=0$ we have strong convergence of $P^{nd}f$ by assumption, and the limit is 
$E(f|\Sigma_{i,d}) = E(f|\Sigma_D)$. Since the atoms $A_0,\dots,A_{d-1}$ generate $\Sigma_D$,
the $\Sigma_D$-measurable functions are of the form $\sum_{\ell=0}^{d-1} c_\ell 1_{A_\ell}$, 
and then 
$E(f|\Sigma_D)=\sum_{\ell=0}^{d-1} \big(\mu(A_\ell)^{-1} \int_{A_\ell} f \,d\mu\big) 1_{A_\ell}$ 
by disjointness of the atoms. This proves \eqref{L1-limits} for $j=0$, since 
$\mu(A_\ell) =d^{-1}$. For $j>0$ we apply $P^j$ and use $P^j 1_{A_\ell}=1_{A_{\ell \dot{+} j}}$.

 Since $P$ is a contraction, \eqref{conditional} follows from 
$\|P^{nd}\big((f-E(f|\Sigma_D)\big)\|_p \to 0$. 
\end{proof}

{\bf Remarks.} 1.  The convergence \eqref{conditional} makes precise \eqref{lly}, giving 
information on the cycle, showing that the permutation of the $y_j$ is cyclic, and expliciting 
the functionals  $\varphi_j$. The assumption of Theorem \ref{periodic-convergence} means 
that $P$ is $L^1$-constrictive, by Corollary \ref{Lly}.

 2. Foguel's assumption in \cite[Theorem 3.3]{F6} yields that both $P^{nj}$ and
$P^{*nj}$ converge strongly in $L^1$ as $n \to \infty$ (his $L^1$-zero-two theorem for $P^{*j}$
 is equivalent to the $L^\infty$-zero-two theorem for $P^j$, and the above convergences 
both hold in the zero alternative). In \cite[Theorem 6.4]{F3} Foguel proved (i) and (iii) assuming
 Harris recurrence. Our result applies without any of these stronger assumptions, also  when
 $P^{*nk}$ does not converge  strongly for any $k$.

3. Wittmann's uniform $L^p$-zero-two theorem \cite[Theorem 1.7]{Wi} yields the strong 
convergence in $L^2$ of $(P^{nj})_{n\ge 1}$  and $(P^{*nj})_{n\ge 1}$ if
$$ 
\lim_{n \to \infty} \|P^n(I-P^j)\|_2 < \sqrt{3}.
$$

4.  A  sufficient spectral condition for the strong convergence in $L^2$ of 
$(P^{nk})_{n\ge 1}$  and $(P^{*nk})_{n\ge 1}$ for {\it some} $k \ge 1$ is given by 
\cite[Proposition 2]{L3}: $\sigma(P) \cap \mathbb T \ne \mathbb T$.
This condition is satisfied if $\lim_{n \to \infty} \|P^n(I-P)\|_2 < 2$ \cite[Theorem 3]{L3}.

5. Derriennic's \cite[Th\'eor\`eme 2]{De}  yields that $\|P^n f\|_2 \to 0$ if and only if
$ \int fg\,d\mu=0$ for every $g \in \bigcap_{n \ge 1} P^{*n}\{h \in L^2: \|h\| \le 1\}$.

\begin{cor} \label{harris}
Let $P$ be an ergodic Harris recurrent Markov operator with invariant probability $\mu$
and period $d$. Then the unimodular eigenvalues of $P$ in $L^1(\mu)$ are precisely
all {\rm d}th roots of unity, the corresponding eigenspaces are one-dimensional,
and the corresponding eigenfunctions have constant absolute value.
\end{cor}
\begin{proof} By \cite{F3} \eqref{L1-limits} holds, so we can apply Theorem 
\ref{periodic-convergence}(ii).
\end{proof}

{\bf Remarks.} 1. \v Sid\'ak \cite{S} proved the result for $P$ on a countable state space
defined by a positive recurrent irreducible Markov matrix. We have not found a reference 
for the corollary in the general case.

2. Gerlach \cite[Theorem 3.4]{Ge}  proved that if $T$ is an irreducible 
power-bounded Harris-type positive operator on a Banach  lattice with order 
continuous norm, then some power $T^n$ has no unimodular eigenvalues different 
from one. This abstract result yields that if $T$ has unimodular eigenvalues, 
they are all $n$th roots of unity, but does not yield the result (which depends on
the cyclically moving sets) that {\it all} $n$th roots of unity are eigenvalues of $T$. 
This follows from Schaefer \cite{Sc2}, since the corresponding eigenfunctions
are bounded by ergodicity.
\medskip
 
It is known that $P$ and $P^*$
 have the same invariant sets, hence the same integrable invariant functions (see 
\cite[Chapter VII]{F2}, for example). In particular, $\Sigma_{i,k}(P)=\Sigma_{i,k}(P^*)$,
and $P^*$ is ergodic when $P$ is. In general $\Sigma_D(P) \ne \Sigma_D(P^*)$,
but $\Sigma_U:=\Sigma_D(P) \bigcap \Sigma_D(P^*)$ generates the unitary subspace
$\mathcal K:=\{f \in L^2: \|P^nf\|_2=\|P^{*n}f\|_2=\|f\|_2 \text{ for every } n \in \mathbb N\}$,
i.e. $\mathcal K=L^2(S,\Sigma_U,\mu)$ \cite[Chapter VIII]{F2}.


\begin{cor} \label{unitary}
Let $P$ be an ergodic bi-stochastic Markov operator as in Theorem  \ref{periodic-convergence},
and let $A_0,\dots,A_{d-1}$ be the atoms generating $\Sigma_D(P)$. Then:

(i)  $\Sigma_U=\Sigma_D(P) =\Sigma_{i,d}$, and $P^* 1_{A_j}=1_{A_{j\dot{-}1}}$.

(ii) $\ P^{*n} \Big(f - d\sum_{\ell=0}^{d-1} \big( \int_{A_\ell}f\,d\mu\big) 1_{A_\ell}\Big) \to 0$ 
 weakly, for every $f \in L^p$, $1\le p< \infty$.
\end{cor}
\begin{proof} By   Theorem \ref{periodic-convergence} and Lemma \ref{k-invariants},
$\Sigma_D(P)=\Sigma_{i,d}(P)=\Sigma_{i,d}(P^*) \subset \Sigma_D(P^*)$,
so $\Sigma_U=\Sigma_D(P)$.  For $f$ in the unitary space $\mathcal K=L^2(\Sigma_U,\mu)$ 
we have  $P^*Pf=f$, hence $P^*1_{A_{j+1}}=P^*P1_{A_j}=1_{A_j}$.
(ii) follows from \eqref{conditional}.
\end{proof}

\begin{theo}
Let $P$ be an ergodic bi-stochastic Markov operator on $(S,\Sigma,\mu)$. Then
$(P^{nd})_{n\ge 1}$ converges in the strong operator topology (SOT) of $L^2$ for some $d$ if and 
only if $P^{*k}P^k$ converges (SOT) to a projection on a finite-dimensional subspace of $L^2$.
\end{theo}
\begin{proof}	
Assume that $P^{nd}$ converges strongly as $n \to \infty$. By Theorem \ref{periodic-convergence},
\eqref{conditional} yields that $\|P^nf\|_2 \to 0$ for every $f \in L^2$ with $E(f|\Sigma_D)=0$.
Hence by \cite[Lemma 3, p. 108]{R}, $P^{*k}P^k$ converges strongly to the projection
$E(f|\Sigma_D)$, which has finite-dimensional range since $\Sigma_D$ is finite by Theorem 
\ref{periodic-convergence}.

Assume now that $P^{*k}P^k$ converges strongly to a projection $E_0$ with finite-dimensional
range. By \cite[Lemma 3, p. 108]{R}, the limit is $E_0f=E(f|\Sigma_D(P))$, and $\|P^nf\|_2 \to 0$ 
for every $f \in L^2$ with $E(f|\Sigma_D(P))=0$, i.e. for every $f \perp E_0L^2$. Since 
$E_0$ has finite-dimensional range, $\Sigma_D(P)$ is finite, therefore atomic. 
Hence the unitary $\sigma$-algebra $\Sigma_U=\Sigma_D(P) \cap \Sigma_D(P^*)$ is finite, 
with $d$ atoms.  By Foguel's \cite[Corollary 2.9]{F6}, $P^{nd}f$ converges {\it weakly}. Since 
$E_0L^2=L^2(S,\Sigma_D(P))$ is finite-dimensional and $P$-invariant, for $g \in E_0L^2$ 
we have that $P^{nd}g$ converges strongly. Together with the above convergence on 
$(E_0L^2)^{\perp}$, we conclude that $P^{nd}$ converges in the strong operator topology 
as $n \to \infty$.
\end{proof}

{\bf Remark.} The necessary condition for convergence of $(P^{nk})_{n\ge 1}$ for some $k$,
that $\Sigma_D$ be finite (see Propositions \ref{foguel} and \ref{deterministic}), implies, by 
\cite[Corollary 2.9]{F6}, that $P^{nd}$ converges in the {\it weak} operator topology (WOT) 
for some $d$.  However, this condition is not sufficient for the SOT convergence (e.g. 
\cite[p. 113]{R}), and is not necessary for WOT convergence (e.g. an invertible mixing 
measure preserving transformation).

\begin{theo} \label{cycle}
Let $P$ be an ergodic $L^2$-quasi-compact bi-stochastic Markov operator on 
$(S,\Sigma,\mu)$ with period $d>1$. Then:

(i) The $\sigma$-algebra  $\Sigma_{i,d}$ of $P^d$-invariant sets is finite,  
has an atom $A_0$ such that for $1\le j\le d-1$ there are disjoint atoms 
$A_j \in \Sigma_{i,d}$ with $P^j1_{A_0}=1_{A_j}$, $P1_{A_j}=1_{A_{j+1}}$ ($A_d=A_0$),
$P^*1_{A_{j+1}}=1_{A_j}$, and $\Sigma_{i,d}$ is generated by $\{A_0,A_1,\dots,A_{d-1}\}$. 

(ii) For every $f \in L^p(\mu)$, $1 \le p< \infty$, and $0 \le j \le d-1$ we have
\begin{equation} \label{Lp-limits}
\lim_{n\to\infty}\big\|P^{nd+j}f - 
d \sum_{\ell=0}^{d-1} \big( \int_{A_\ell} f \,d\mu\big) 1_{A_\ell\dot{+}j} \big\|_p =0,
\end{equation}

\begin{equation} \label{P*-limits}
\lim_{n\to\infty}\big\|(P^*)^{nd+j}f - 
d \sum_{\ell=0}^{d-1} \big( \int_{A_\ell} f \,d\mu\big) 1_{A_\ell\dot{-}j} \big\|_p =0.
\end{equation}

(iii) The unimodular eigenvalues of $P$ on $L^p$ are precisely all {\rm d}th roots of unity, 
and the corresponding eigenfunctions have constant absolute value.
\end{theo}
\begin{proof}
By Proposition \ref{qc}, we have $L^2$ operator norm convergence of $(P^{nd})_{n\ge 1}$,
and since $P^*$ is clearly $L^2$-quasi-compact with the same period, also $(P^{*nd})_{n\ge 1}$
converges in operator norm.

Let $d_0$ be the smallest $k$ such that $(P^{nk})_{n \ge 1}$ converges strongly in $L^2$.
By Proposition \ref{divisor}, $d_0|d$. We show $d_0=d$: observe that for every unimodular
 eigenvalue $\lambda$ we have, by definition, convergence  $\lambda^{nd_0}$. Since all
these eigenvalues are roots of $d$, the minimality of $d$ implies $d|d_0$.

We now obtain (i) and \eqref{Lp-limits} from Theorem \ref{periodic-convergence} and Corollary
\ref{unitary}; \eqref{P*-limits} follows by applying \ref{periodic-convergence} to $P^*$.

By the construction of Corollary \ref{ue} and the definition of $d$, all unimodular eigenvalues 
of $P$ are $d$th roots of unity. (iii) follows from Theorem \ref{periodic-convergence}(ii).
\end{proof}

We now extend to the periodic case the limit theorem for $P^n$, proved in Proposition \ref{exponential}(ii) for the aperiodic case.
\begin{theo}
 Let $P$ be an ergodic $L^2$-quasi-compact bi-stochastic Markov operator with period $d>1$, 
and let $A_0,\dots,A_{d-1}$ be the atoms generating $\Sigma_{i,d}$. 

(i) If $f \in L^r(\mu)$, $1 \le r < \infty$, with $E(f|\Sigma_{i,d})=0$ (i.e. 
$\int_{A_j}f\,d\mu=0$ for $0\le j\le d-1$), then $\|P^nf\|_r \to 0$ and $\|P^{*n}f\|_r \to 0$.

(ii) If $r>1$ and $f \in L^r(\mu)$ with $E(f|\Sigma_{i,d})=0$, then $P^nf \to 0$ and 
$P^{*n}f \to 0$ a.e.

(iii) Assume $P$ is hyperbounded and maps $L^1(\mu)$ to $L^q(\mu)$ for some $q>1$. 
If $f \in L^1(\mu)$ satisfies $E(f|\Sigma_{i,d})=0$, then $P^nf \to 0$ a.e.
\end{theo}
\begin{proof}
Since $P^d$ and $P^{*d}$ have the same invariant sets, it is enough to prove the assertions
for $P$, and then apply them to $P^*$ and obtain the convergence results for powers of $P^*$.
\smallskip

(i) First assume $r>1$. The assumption on $f$ means $E_df=0$, so by Proposition 
\ref{exponential}(i) we have $\|P^{nd}f\|_r \to 0$.  Since $P$ is a contraction, 
$\|P^{n}f\|_r \to 0$.  

For $f$ bounded with $E_df=0$, we have  $\|P^nf\|_1 \le \|P^nf\|_r \to 0$.
Standard approximations (using $\|E_d\|_\infty=1$) yield the convergence for $f \in L^1$ 
with $E_df=0$.

(ii) For $1<r< \infty$, Proposition \ref{exponential}(i) yields 
$\|P^{dn}f\|_r \le C_r\rho_r^n \|f\|_r$ when $E_df=0$, so
$\sum_{n=1}^\infty \|P^{nd}f\|_1 < \infty$.  Since $P$ is a contraction,
$\sum_{n=1}^\infty \|P^{nd +j}f\|_1 < \infty$ for $0\le j \le d-1$. Hence
$$
\sum_{n=0}^\infty \|P^{n}f\|_1 =
\sum_{j=0}^{d-1}\sum_{n=0}^\infty \|P^{nd +j}f\|_1 < \infty.
 $$
By Beppo Levi $\sum_{n=0}^\infty |P^{n}f| < \infty$ a.e., so $P^{n}f \to 0$ a.e. 

(iii) The proof is similar to the proof of (ii), using Corollary \ref{doeblin}(i).
\end{proof}

{\bf Remarks.} 1. In case (iii), $P^*$ need not map $L^1$ to some $L^r$, although $P$ does, 
so for $L^1$ functions we obtain from the theorem only a.e. convergence of $P^nf$, but not of
 $P^{*n}f$. However, as noted in the  remarks following Example 1, in case (iii) $P$ is Harris, 
so on each invariant set of $P^d$ the restriction of $P^d$ and $P^{*d}$ are ergodic and Harris, 
and by Horowitz \cite{H}, with the assumptions $ f \in L^1$ and $E_df=0$,  we obtain 
$P^nf \to 0$ a.e. and $P^{*n}f \to 0$ a.e.

2. The theorem applies to $P$ ergodic hyperbounded.
In general,  a hyperbounded $P$ need not be  Harris (see Theorem \ref{riesz} below), 
so only (i) and (ii) of the theorem apply.

3. Equivalence of finiteness of $\Sigma_D$ and  part (i) of the theorem was proved in 
\cite{BK} for kernel operators as in Example 1 (which are clearly Harris recurrent), 
without any further assumptions on the kernel $k(x,y)$, i.e. without assuming hyperboundedness. 
Additional equivalent conditions are also given there.
\medskip

\noindent
{\bf Problem 1.} {\it If $P$ is ergodic hyperbounded, does $P^nf \to 0$ a.e for $f \in L^1$ 
with $E_df=0$?}

For  the aperiodic case, the problem is whether $P^nf$ converges a.e. for every $f \in L^1$.
Note that for aperiodic ergodic Harris recurrent operators, Horowitz \cite{H} proved a.e. 
convergence of $P^n f$ for every $f \in L^1$.
\medskip

\section{Conditions for aperiodicity of hyperbounded Markov operators}

In this section we look for conditions of aperiodicity of an ergodic hyperbounded (necessarily
bi-stochastic) Markov operator. By the definition, when $d=1$ we have 
$\sigma(P) \cap \mathbb T=\{1\}$ (for $P$ on $L^p$ which  maps $L^p$ to $L^q$), and 
 with the uniform ergodicity in Corollary \ref{ue}, $\|P^n-E\|_p \to 0$  by \cite[Theorem 4]{Lu}. 
Obviously, if $P^n$ converges in operator norm, the hyperbounded operator is aperiodic. 
In the complex $L^p$, this convergence is equivalent to a (global) spectral gap:
$r(P_{|(I-P)L^p})< 1$. For any ergodic Markov operator preserving $\mu$,
$\overline{(I-P)L^p} = \{f \in L^p: \int f\,d\mu=0\}$ by the mean ergodic theorem.

Recall that by Proposition \ref{all-Lp}, a hyperbounded Markov operator is hyperbounded in 
each $L^r$, $1<r< \infty$, and the period is the same in all $L^r$, by Proposition \ref{qc-Lp}.
We therefore look at hyperboundedness in $L^2$.

\begin{theo} \label{wang}
Let $P$ be an ergodic hyperbounded Markov operator mapping $L^2(\mu)$ to $L^4(\mu)$.
If $\|P\|_{L^2 \to L^4} < 2^{1/4}$, then $P$ is aperiodic.
\end{theo}
\begin{proof}
We note that for real functions it is easy to prove (since $P$ is real), that
\begin{equation} \label{real-gap}
\sup\{\|Pf\|_2: \|f\|_2=1, \int f\,d\mu=0\} =
\end{equation}
$$
\sup\{\|P(f+ig)\|_2: \|f+ig\|_2=1, \int (f+ig) d\mu=0\} .
$$
Since eigenfunctions of unimodular eigenvalues different from 1 have integral zero, to prove 
aperiodicity it is enough to prove that the left-hand side of \eqref{real-gap} is less than 1.

The proof is the same as Wang's proof of \cite[Theorem 1.1]{W1}. Wang assumes that $P$
is symmetric, but this is because he assumes $P$ to contract only $L^2(\mu)$, with $P1=1$;
symmetry is used in his proof only for obtaining $P^*1=1$, which holds for any bi-stochastic
Markov operator $P$ since $\mu$ is an invariant probability.
\end{proof}

\begin{prop} \label{optimal}
There exists an ergodic hyperbounded Markov operator $P$ with period $d >1$,
which for every $q>2$ maps $L^2(\mu)$ to $L^q(\mu)$, 
with $\|P\|_{L^2 \to L^q} = 2^{\frac12-\frac1q}$. 
\end{prop} 
\begin{proof} Let $P$ be the Markov operator of Example 2, which is ergodic and symmetric
with period 2. The explicit definition of $P$, by Example 1, is
$Pf = 2a1_{[0,1/2)} + 2b1_{[1/2,1]}$, where $a=\int_{1/2}^1 f(y)dy$ and
$b= \int_0^{1/2} f(y)dy$. Then (for real $f$)
$\|Pf\|_2^2=4a^2\cdot\frac12 +4b^2\cdot\frac12= 2(a^2+b^2)$.
Hence for $2<q<\infty$ we have
$$
\|Pf\|_q =\big(2^q |a|^q\cdot\frac12 +2^q |b|^q\cdot\frac12\big)^{1/q}=
2\cdot 2^{-1/q} (a^q +b^q)^{1/q} \le 
$$
$$
2^{1-1/q} (a^2 + b^2)^{1/2} =
2^{1-1/q} 2^{-1/2} \|Pf\|_2 \le 2^{\frac12 -\frac1q} \|f\|_2.
$$
Hence   $\|P\|_{L^2 \to L^q} \le 2^{\frac12-\frac1q}$.

Define $f=2\cdot 1_{[0,1/2)}$. Then $\|f\|_2= \sqrt 2$, and 
$\|Pf\|_q= 2^{1-\frac1q} = 2^{\frac12-\frac1q}\|f\|_2$. Thus, for $q< \infty$,
$\|P\|_{L^2 \to L^q} = 2^{\frac12-\frac1q}$.

Now let $q=\infty$. Then, for any $f \in L^2$,
$$
\|Pf\|_\infty= 2\max\{|a|,|b|\} \le 2 \sqrt{a^2+b^2}=\sqrt 2 \|f\|_2,
$$
which shows $\|P\|_{L^2 \to L^\infty} \le \sqrt 2$. For $f =2\cdot 1_{[0,1/2)}$ we
have $\|Pf\|_\infty = 2 =\sqrt 2 \|f\|_2$; hence $\|P\|_{L^2 \to L^\infty}= \sqrt 2$.
\end{proof}

{\bf Remarks.} 1. For $q=4$ the proposition shows the optimality of the constant $2^{1/4}$
in Theorem \ref{wang} (and in \cite[Theorem 1.1]{W1}). Wang's simple example 
\cite[p. 2633]{W1} is not ergodic. 
Wang's definition of spectral gap, the norm on functions of integral 0 being less than 1, 
is valid only for the ergodic case. In general, the norm should be taken on the 
subspace orthogonal to the invariant functions, and in this sense Wang's example has a 
spectral gap ($P^n$ converges in norm).

2. Combining Theorem \ref{wang} with Proposition \ref{optimal}, we obtain that for any
$q>4$ there exists $\Delta_q \in[2^{\frac14},2^{\frac12 -\frac1q}]$ such that if $P$
is an ergodic Markov operator mapping $L^2$ to $L^q$ with $\|P\|_{L^2 \to L^q} < \Delta_q$,
then $P$ is aperiodic.

3. Let $q \in (2,4)$. If there exists $\Delta_q$, such that every ergodic $P$ mapping
$L^2$ to $L^q$ with $\|P\|_{L^2 \to L^q} < \Delta_q$ is aperiodic, then 
$ \Delta_q \le 2^{\frac12 -\frac1q}$,  by Proposition \ref{optimal}.

\begin{theo} \label{wang3}
Let $P$ be an ergodic hyperbounded Markov operator mapping $L^2(\mu)$ to $L^3(\mu)$.
If $\|P\|_{L^2 \to L^3} < 2^{1/6}$, then $P$ is aperiodic. The value $2^{1/6}$ is optimal.
\end{theo}
\begin{proof}
We prove that the left-hand side of \eqref{real-gap}, denoted later by $\rho$,
 is less than 1. Thus, we consider real functions.  Our proof is inspired by \cite{W1}.

 Fix $f \in L^2(\mu)$ with $\|f\|_2=1$ and $ \mu(f):= \int f d\mu=0$. 
We may assume $\mu((Pf)^3) \ge 0$ (otherwise we replace $f$ by $-f$).

Let $\varepsilon \in (0,1)$, and define $g = \varepsilon^{1/2} +(1-\varepsilon)^{1/2} f$. 
Then  $\|g\|_2^2=  \varepsilon +(1-\varepsilon)\|f\|_2^2 = 1 $, by orthogonality. Then
$$
\delta:= \|P\|_{L^2 \to L^3}^3 \ge  \|Pg\|_3^3 \ge \int (Pg)^3d\mu= 
$$
$$
 \int  \Big(\varepsilon^{3/2} +3 \varepsilon(1-\varepsilon)^{1/2}Pf +
3 \varepsilon^{1/2}(1-\varepsilon)(Pf)^2 +(1-\varepsilon)^{3/2}(Pf)^3 \Big) d\mu \ge
$$
$$
\varepsilon^{3/2} +3\varepsilon^{1/2}(1-\varepsilon) \mu((Pf)^2) .
$$
Hence $\displaystyle{\|Pf\|_2^2 \le 
\frac{\delta-  \varepsilon^{3/2}}{3 \varepsilon^{1/2}(1-\varepsilon) } }$.
Since $f$ with norm 1 and zero integral was arbitrary, we obtain $\displaystyle{\rho^2 \le
\frac{\delta-  \varepsilon^{3/2}}{3 \varepsilon^{1/2}(1-\varepsilon) } }$.
Taking $\varepsilon=\frac12$ we find 
$\rho^2 \le (\delta- (\frac12)^{3/2})/ 3(\frac12)^{1/2} \cdot \frac12  $.  We have  
$$ 
\frac{\delta- (\frac12)^{3/2}}{3(\frac12)^{3/2} } <1 \ \text{if and only if }
\delta < 4\cdot (\frac12)^{3/2} = \sqrt 2. 
$$ 
Hence $\|P\|_{L^2 \to L^3} < 2^{1/6}$ implies $\rho^2 <1$, which yields aperiodicity.
\smallskip

The optimality of the value $2^{1/6}$ follows from the example in Proposition \ref{optimal},
with $q=3$.
\end{proof}

{\bf Remark.} Theorem \ref{wang3} does not follow from Theorem \ref{wang}, since $P$
may map $L^2$ to $L^3$ and not to $L^4$. 
\medskip

\noindent
{\bf Problem 2.} {\it For $q \in (2,3)$, is there a constant $\Delta_q$ such that any ergodic
 hyperbounded $P$ mapping $L^2$ to $L^q$ with $\|P\|_{L^2 \to L^q} <\Delta_q$ is aperiodic?}
Since for $P$ mapping $L^2$ to $L^q$ ($q>2$) we have $\|P\|_{L^2 \to L^q} \ge 1$,
if $\Delta_q$ exists, then $\Delta_q >1$. If $\Delta_q$ does not exist, we may still ask: 
{\it Is any hyperbounded $P$ mapping $L^2$ to $L^q$ with $\|P\|_{L^2 \to L^q} = 1$ aperiodic?}
\medskip

\section{Hyperbounded Markov operators defined by convolutions}

An important class of Markov operators with an invariant probability is given by convolution 
operators on compact groups. We deal in this section with the circle group $\mathbb T$,
identified with $[-\pi,\pi)$, with normalized Haar (Lebesgue) measure $\mu$. 
We use the notation $L^p$ or $L^p(\mathbb T)$  for $L^p(\mathbb T,\mu)$.

Let $\nu$ be a probability on $\mathbb T$, and define the convolution operator
$P_\nu f=\nu*f$, with $\nu*f(x):= \int_\mathbb T f(x-y) d\nu(y)$.  Then $P_\nu1=1$, and for
$f \in L^p$, $1 \le p < \infty$ we have
$$
\|P_\nu f\|_p^p =\int \Big|\int f(x-y)d\nu(y)\Big|^p d\mu(x) \le
\int \int |f(x-y)|^p d\nu(y) d\mu(x) =
$$
$$
=\int \Big(\int |f(x-y)|^p d\mu(x) \Big) d\nu(y) =
\int \int |f(z)|^p d\mu(z) d\nu(y) = \|f\|_p^p.
$$
Similarly $\int P_\nu f\, d\mu=\int f\,d\mu$, so $P_\nu$ is a Markov operator with $\mu$ invariant.

We denote $e_n(x):= \e^{inx}$, $n \in \mathbb Z$. The Fourier coefficients of $f \in L^1$
are $\hat f(n) := \frac1{2\pi} \int f(x)e_{-n}(x)dx$,
and the Fourier-Stieltjes coefficients of $\nu$ are $\hat \nu(n) = \int e_{-n}(x)d\nu(x)$. 
Note that if $\nu << \mu$ with $\frac{d\nu}{d\mu}= \phi \in L^1$, then $\hat \nu =\hat \phi$.
Using Fubini's theorem, we obtain
$$
\widehat{P_\nu f }(n)=\int \int f(x-y)e_{-n}(x)d\mu(x)d\nu(y)= 
\int \int f(z)e_n(z+y)d\mu(z)d\nu(y) =\hat\nu(n)\cdot \hat f(n).
$$  
This yields that $\{\hat\nu(n)\}_{n \in \mathbb Z}$ is a {\it multiplier} sequence in any $L^p$,
i.e., if $f =\sum_{n \in \mathbb Z} c_ne_n \in L^p$,  $1<p< \infty$, then 
$\sum_{n \in \mathbb Z}\hat\nu(n)c_ne_n  \in L^p$ (convergence in $L^p$). 

Hyperboundedness of $P_\nu$, mapping $L^p$ into $L^q$ with $p<q$, means that
 $\{\hat\nu(n)\}_{n \in \mathbb Z}$ in an $L^p-L^q$ multiplier sequence. In some harmonic
analysis papers (e.g. \cite{Rit}, \cite{GHR}, \cite{DHR}), when $P_\nu$ is hyperbounded
$\nu$  is called {\it $L^p$-improving}.

If $P_\nu$ is hyperbbounded in $L^p(\mathbb T)$ for $1\le p< \infty$, then it is hyperbounded 
in every $L^r$, $1 < r< \infty$ (Proposition \ref{all-Lp}), and  quasi-compact in each $L^r$ by
Corollary \ref{ue}. 
Obviously, for any probability $\nu$, the spectrum of of $P_\nu$ in
$L^2$ is $\overline{\{\hat\nu(n): n\in \mathbb Z\}}$. Graham, Hare and Ritter 
\cite[Theorem 4.1]{GHR} proved that {\it if $P_\nu$ is hyperbounded, then for any $1<r< \infty$
the spectrum of $P_\nu$ on $L^r$ is  $\overline{\{\hat\nu(n): n\in \mathbb Z\}}$}.
\medskip

It was shown in \cite[Corollary 4.2(iv)]{DL} that $P_\nu$ is uniformly ergodic in $L^2(\mu)$
if and only if $\inf_{n \ne 0} |\hat\nu(n)-1| >0$. It was proved  \cite[Theorem 4.6]{DL} 
that if $\nu$ is an adapted discrete probability on $\mathbb T$, then $P_\nu$ is not uniformly 
ergodic in $L^2$.  An example of $\nu$ continuous singular satisfying the above  condition,
with $P_\nu$ not Harris, was presented in \cite[Proposition 4.7]{DL}. The discussion in
\cite[p. 92]{DL} shows the existence of a continuous probability with all its powers singular, 
with $P_\nu$ not uniformly ergodic in $L^2$. By \cite[Theorem 4.3]{DL},
if $P_\nu$ is Harris recurrent (some power $\nu^k$ is not singular), then it is uniformly
ergodic in $L^2$. 

A Rajchman probability $\nu$ (i.e. $\hat\nu(n) \to 0$ as $|n| \to \infty$; 
necessarily continuous by a result of Wiener \cite [Theorem III.9.6]{Z})  defines $P_\nu$
uniformly ergodic on $L^2$, by \cite[Theorem 4.4]{DL}.
 The singular probability $\nu$ with $P_\nu$ uniformly ergodic and not Harris, 
presented in \cite{DL}, {\it is not} Rajchman.

\begin{prop} \label{ritter}
Let $\nu$ be a probability on $\mathbb T$. If $(P_\nu)^m$ is hyperbonded for some $m>1$,
then $P_\nu$ is hyperbounded.
\end{prop}
\begin{proof} We first note that $(P_\nu)^m =P_{\nu^m}$. By Corollary \ref{p-2},
 $P_{\nu^m}$ maps $L^p$ to $L^2$, for some $1<p<2$. By Ritter \cite[Lemma 2]{Rit},
$P_\nu$ is hyperbounded.
\end{proof}

{\bf Remark.} As noted in the remarks preceding Proposition \ref{fixed-space}, Proposition
\ref{ritter} does not extend to general bi-stochastic Markov operators.

\begin{theo} \label{riesz}
There exists a singular Rajchman probability $\nu$ on $\mathbb T$ such that all its 
convolution powers are singular (so $P_\nu$ is not Harris) and $P_\nu$ is hyperbounded. 
\end{theo} 
\begin{proof} We use the construction of  Riesz products (see \cite[Section V.7]{Z}),
used by Zafran \cite[p. 619]{Za}.  Let $-1\le a_k \le 1$ with $a_k\not=0$ for every 
$k\ge 1$. Let $(n_k)$ be a lacunary sequence, with $n_{k+1}/n_k\ge q>3$, for every 
$k\ge 1$. Expanding the partial Riesz products, with $m_n =\sum_{k=1}^n n_k$, we obtain
$$
p_n(x)=\Pi_{k=1}^n(1+a_k\cos(2\pi n_k x)) = 1+ \sum_{j=1}^{m_n} \gamma_j \cos(2\pi jx),
$$
which is a partial sum of the Fourier-Stieltjes series 
$1+\sum_{j=1}^\infty \gamma_j\cos(2\pi jx)$ of a non-decreasing continuous function
 $F(x)$ \cite[Theorem V.7.5, p. 209]{Z}, given by the relation
$$
F(x)-F(0)=\lim_{n\to\infty}\int_{0}^x  p_n(t)dt.
$$
Since $p_n(t) \ge 0$ and $\int_0^1p_n(t)=1$, we have $F(1)-F(0)=1$.
The function $F$ has zero derivative Lebesgue a.e., i.e. $F$  is singular,
if $\sum_{k=1}^\infty a_k^2=\infty$; otherwise, $F$ is absolutely continuous
(\cite[Theorems  V.7.6 and  V.7.12]{Z}).

By expanding the partial Riesz  products, it is simple to see that
$\gamma_m=0$ if $m$ is not of the form $\sum_{j=1}^k \epsilon_j n_j$,
where $\epsilon_j \in \{-1,0,1\}$. Let $\nu$ be the positive measure defined by $F$.
Since $F(1)-F(0)=1$, $\nu$ is a probability, and  $\hat\nu(m) = \gamma_m$. 
By \cite{Za}, we have
$$
\hat \nu(\sum_{j=1}^k \epsilon_j n_j )=\Pi_{j=1}^k (\frac{a_j}{2})^{|\epsilon_j|}.
$$
We conclude that if $\nu$ is Rajchman, we must have $a_k =2 \hat\nu(n_k) \to 0$,
while if $a_k\to 0$, then $\nu$ is Rajchman. Since $|a_k| \le 1$ we have 
 $|\hat\nu(m)| \le \frac1{2}$ for $m \ne 0$, so by the theorem of Ritter for Riesz products 
\cite{Rit}, $P_\nu$ is hyperbounded.

For $n \ge 1$, by \cite[p. 619]{Za} the $n$-fold convolution power $\nu^n $ 
is represented by the infinite Riesz product:
$$
\Pi_{k=1}^\infty(1+2(\frac{a_k}2)^n\cos(2\pi n_k x)).
$$
By \cite{Z}, $\nu^n$ is singular if and only if
$ \sum_{k=1}^\infty a_k^{2n} = \infty$. By choosing $(a_k)$ such that
$\sum_{k=1}^\infty |a_k|^n =\infty$ for every $n \ge 1$, e.g. $a_k =1/\log(k+2)$
for $k \ge 1$, we make sure that all convolution powers of $\nu$ are singular,
and that $\nu$ is Rajchman. 
\end{proof}

{\bf Remark.} The existence of a singular Rajchman probability with all its powers singular
 is a special case of Theorem S of Varopoulos \cite{V}.  Our proof for $\mathbb T$ 
is along classical lines, using Riesz products  and Zygmund's criteria for singularity, 
and  yields a concrete $\nu$; it is simpler than the proof for general LCA groups in \cite{V}, 
which  shows only the {\it existence} of the desired measures.

\begin{cor} \label{not-hyper}
There exists a singular Rajchman probability $\nu$  such that all its convolution
 powers are singular and $P_\nu$ is not hyperbounded.
\end{cor}
\begin{proof} We denote by $\nu_0$ the Rajchman  probability given by Theorem
 \ref{riesz}. Let $(a_n)_{n\ge 0}$ be a sequence of positive numbers with $\sup_n\, a_n<1$
 and $a_n \to 0$. Badea and M\"uller \cite[Theorem 5]{BM} proved that there exists
a probability measure $\nu << \nu_0$ with $|\hat\nu(n)| \ge a_{|n|}$ for every integer
$n$. By a result of Rajchman (see references and a proof in \cite[(2.1)]{Ly2}), also $\nu$ is
Rajchman. Obviously $\nu^k<<\nu_0^k$, so  $\nu^k$ is singular. In order to show that we
can obtain $P_\nu$ not hyperbounded, we choose $a_n \to  0$ very slowly, so that Edwards's
necessary conditions \cite[formulas (16.4.9) and (16.4.10)]{Ed2} are violated; e.g.
$a_n =1/\log\log(n+27)$ for $n \ge 0$.
\end{proof}

{\bf Remarks.} 1. It is noted in \cite[p. 487]{GHR} that it is possible to construct a Rajchman
probability $\nu_0$ such that for {\it  any} probability $\nu << \nu_0$ the operator $P_\nu$
 is not hyperbounded. Our construction is different.

2. Since $\nu$ of the corollary is Rajchman, $P_\nu^nf= \nu^n*f \to \int f\,d\mu$ a.e for every 
$f \in L^p$, $p>1$, \cite{DL}, although $P_\nu$ is not hyperbounded and not Harris recurrent.

3.  Sarnak \cite[p. 309]{Sa} proved that {\it if $\nu$ is Rajchman, then for any $1<r< \infty$
the spectrum of $P_\nu$ on $L^r$ is  $\overline{\{\hat\nu(n): n\in \mathbb Z\}}$},
even when, as in Corollary \ref{not-hyper}, \cite[Theorem 4.1]{GHR} does not apply.

\begin{prop} \label{Lr-deriv}
Let $\nu$ be a probability on $\mathbb T$ and $1<q<\infty$. The operator $P_\nu$ maps $L^1$ into
$L^q$ if and only if $\nu<<\mu$ with $\frac{d\nu}{d\mu} \in L^q$.
\end{prop}
\begin{proof} Let $\nu<<\mu$ with $\phi=\frac{d\nu}{d\mu} \in L^q$. 
Putting $p=1$ and $r=q$ in \cite[Theorem II.1.15]{Z}, we obtain that for
 $f \in L^1$ we have $\phi*f \in L^q$, and
$$
\|P_\nu f\|_q =\|\phi*f\|_q \le \|\phi\|_q\cdot \|f\|_1.
$$
Hence $P_\nu$ maps $L^1$ into $L^q$, with norm $\|P_\nu\|_{L^1\to L^q} \le \|\phi\|_q$.
\smallskip

The converse is Theorem 16.3.4 of \cite{Ed2}.
\end{proof}

{\bf Remarks.} 1. When $\nu<<\mu$,  $P_\nu$ is an integral operator, defined by the kernel
$k(x,y)=\phi(x-y)$. Similarly, if $\nu$ is not singular (with respect to $\mu$), then
$P_\nu$ bounds a non-zero integral operator, and therefore is Harris recurrent.

2. The probability $\nu$ of  Theorem \ref{riesz} yields $P_\nu$ which is hyperbounded in
every $L^p$, $1<p< \infty$, but not in $L^1$ since $\nu$ is singular.

\begin{prop} \label{ghr}
Let $\nu$ be a probability on $\mathbb T$, such that $P_\nu$ is  hyperbounded, mapping 
$L^p$ to $L^q$, $1<p<	q$, and put $\alpha =q/p$. If $r>\alpha/(\alpha-1)$ and
 $\eta<<\nu$ is a probability with $\frac{d\eta}{d\nu} \in L^r(\nu)$, then $P_\eta$ is
hyperbounded.
\end{prop}
\begin{proof} Let $s=r/(r-1)$. Then $1<s< \alpha$, and the proof of \cite[Theorem 1.1]{GHR}
shows that $P_\nu$ maps $L^{ps}$ to $L^q$.
\end{proof}

\begin{prop}
There exists a singular probability $\nu$ on $\mathbb T$ which is {\rm not} Rajchman, 
such that all its convolution powers are singular, and $P_\nu$ is hyperbounded. 
\end{prop} 
\begin{proof}
Let $\nu$ be the usual Cantor-Lebesgue measure, as described in \cite[Proposition 4.7]{DL}.
It is clearly not Rajchman, since $\hat\nu(n)$ is constant along the powers of 3, and it is
shown there that all powers of $\nu$ are singular. Oberlin \cite{Ob} proved that $\nu$ is
 hyperbounded ($L^p$-improving); see also \cite[Proposition 4.2]{DHR}.
\end{proof}

{\bf Remarks.} 1. It was proved by Graham et al. \cite[Corollary 3.2(iii)]{GHR} that if $\nu$ maps
{\it every} $L^p$, $1<p<2$, into $L^2$, then $\nu$ is Rajchman.

2. If $P_\nu$ is hyperbounded, then $\limsup_{|n| \to \infty} |\hat\nu(n)| <1$: by Proposition
\ref{all-Lp} and duality $P_\nu$ maps some $L^p$, $1<p<2$, to $L^2$, and then we apply
\cite[Corollary 3.2(ii)]{GHR}.

\begin{prop} \label{continuous}
Let $\nu$ be a probability on $\mathbb T$. If $P_\nu$ is hyperbounded, then $\nu$ is 
continuous (atomless), and
\begin{equation} \label{HT}
\sum_{|n|\ne 0} \frac{|\hat\nu(n)|^2}{|n|} < \infty.
\end{equation}
\end{prop}
\begin{proof}
It suffices to prove \eqref{HT}, since then by Kronecker's lemma 
$\frac1{2N+1} \sum_{|n|\le N} |\hat\nu(n)|^2 \to 0$ as  $N\to \infty$, so
by Wiener's  criterion \cite[Theorem III.9.6]{Z}, $\nu$ is atomless.
\smallskip

By Corollary \ref{p-2}, hyperboundedness of $P_\nu$ implies that $P_\nu$ maps some 
$L^p$, $1 < p<2$, into $L^2$. This then allows the following proof, based on an idea of
Hare and Roginskaya \cite{HR}. We show that for any $\varepsilon >0$ we have 
\begin{equation} \label{hr}
\sum_{|n|\ne 0}\frac { |\hat\nu(n)|^2}{n^{2(p-1)/p}\log^{1+\varepsilon}(1+|n|) }< \infty.
\end{equation}
Denote $s=2(p-1)/p$. Let $D_N(x)=\sum_{|n| \le N} \e^{inx}$ be the complex Dirichlet kernel.
By \cite[Lemma 2.1]{AAJRS} we have $\|D_N\|_p^p \le C_p N^{p-1}$ for $p>1$.  Then
$$
\sum_{|n|\ne 0}\frac { |\hat\nu(n)|^2}{n^s \log^{1+\varepsilon}(1+|n|) } =
\sum_{k=0}^\infty
 \sum_{|n|=2^k}^{2^{k+1}}\frac { |\hat\nu(n)|^2}{n^s \log^{1+\varepsilon}(1+|n|) } \le
$$
$$
c_1 +c_2 \sum_{k=1}^\infty \frac1{2^{ks} k^{1+\varepsilon} }\|\nu *D_{2^{k+1}}\|_2^2 \le
  c_1 +c_3 \sum_{k=1}^\infty  \frac1{2^{ks} k^{1+\varepsilon} }2^{2(p-1)/p} =
 c_1 +c_3 \sum_{k=1}^\infty  \frac1{ k^{1+\varepsilon} }< \infty.
$$
This proves \eqref{hr}, and since $s<1$, \eqref{hr} implies \eqref{HT}.
\end{proof}

{\bf Remarks.} 1. Equation \eqref{hr} improves the necessary condition of Edwards 
\cite[Remark (3), p. 302]{Ed2} for $q=2$.

2. For any positive sequence $(a_n)_{n \in \mathbb Z}$ which tends to $0$ as $|n| \to \infty$, 
there exists a  {\it complex}  measure $\eta$ satisfying 
$\sum_{n\in \mathbb Z} a_n|\hat\eta(n)|^2 <\infty$, with $P_\eta f =\eta*f$ not $L^p$ 
improving \cite[Proposition 2.7]{HR}.

3. Another proof of continuity (in a stronger form) is in \cite [Corollary 3.2(i)]{GHR}.

\smallskip

{\bf Problem 3.} {\it Let $(a_n)_{n \in \mathbb Z}$ be a positive sequence which 
tends to $0$ as $|n| \to \infty$. Does there exist a  {\em probability}  measure $\nu$ satisfying 
$\sum_{n\in \mathbb Z} a_n|\hat\nu(n)|^2 <\infty$, with $P_\nu$ not hyperbounded?}
A positive answer for $a_n=1/|n|$ is given below. As mentioned above, a complex measure
with the desired properties exists by \cite[Proposition 2.7]{HR} (and it is possible also to obtain
from it  a real {\it signed} measure). We conjecture that the answer to Problem 3 is always
positive; if not, it means that {\it there exists a positive sequence $a_n \to 0$ as $|n| \to \infty$
(necessarily with $\sum_n a_n=\infty$), such that for every probability $\nu$ with 
$\sum_{n\in \mathbb Z} a_n|\hat\nu(n)|^2 <\infty$, the convolution operator $P_\nu$ is hyperbounded.} Based on the current knowledge, such a sufficient condition for hyperboundedness
seems unlikely.

\begin{prop} \label{not-hyper2}
There exists on $\mathbb T$ a probability $\nu << \mu$ such that $P_\nu$ is not
hyperbounded, but \eqref{HT} holds.
\end{prop}
\begin{proof}
Let $b_n=1/\log(|n|+2)$ for $n \in \mathbb Z$. Then for $n>0$ we have
$b_{n-1}+b_{n+1}-2b_n \ge 0$, so by \cite[Theorem I.4.1]{Ka} (see also 
\cite[Theorem V.1.5]{Z}), there is a {\it non-negative} function $f \in L^1(\mathbb T)$ 
such that $\hat f(n)=b_n$, $n \in \mathbb Z$.
Let $g = f/\|f\|_1$ and define $d\nu=g\,d\mu$. Since $\hat\nu(n) =1/\|f||_1 \log(|n|+2)$,
the series in \eqref{HT} converges.

However, for any $1<p<2$ we have $2(p-1)/p <1$, so the series in \eqref{hr} always diverges,
which implies that $P_\nu$ is not hyperbounded.
\end{proof}

 {\bf Remarks.} 1. Corollary \ref{not-hyper} gives an example of $\nu$ Rajchman singular
 with a non-Harris $P_\nu$ which is  $L^2$ uniformly ergodic and not hyperbounded.
 Proposition \ref{not-hyper2} provides an example with $\nu$ absolutely continuous 
(so $P_\nu$ is Harris).

2. It is easy to find $\nu$ not continuous with $P_\nu$ uniformly ergodic in $L^2$
 and not hyperbounded: take $\nu_1<<\mu$ with $d\nu_1/d\mu$ in $L^2$, and define
 $\nu=\frac12(\delta_0+\nu_1)$. Since $\delta_0 * f=f$ for any $f$, we have $P_\nu=\frac12(I+P_{\nu_1})$. By \cite{DL} $P_{\nu_1}$ is uniformly ergodic in $L^2$ 
(and also hyperbounded by Proposition \ref{Lr-deriv}).   Also $P_\nu$ is
uniformly ergodic (see \cite{DL}),  but is not hyperbounded by Proposition \ref{continuous}.
See also the proof of Proposition \ref{no-hyper}, with $Q=P_{\nu_1}$ and $P=P_\nu\,$.

3. Since $P_\mu$ is hyperbounded, the first part of Proposition \ref{not-hyper2} follows from 
\cite[Theorem 1.4]{GHR}. Another construction is given in \cite[Theorem 2.3]{GHR}.
  Our proof is different. The construction of Corollary \ref{not-hyper} shows that for every
Rajchman probability $\nu_0$ there exists $\nu<<\nu_0$ with $P_\nu$ not hyperbounded,
since \eqref{HT} fails.

\begin{theo} \label{sum-estimates}
Let $\nu$ be a probability on $\mathbb T$, such that $P_\nu$ is hyperbounded,
mapping $L^p$ to $L^q$, with $1 \le p<q < \infty$. 

(i)  If $1 \le p<2$, then we may assume $q\le 2$, and then
there exists a constant $C_{p}$ such that 
for every $a \in \mathbb Z$ and $b \in \mathbb N^+$ we have
\begin{equation} \label{sum-estimate1}
\sum_{|n| \le N}|\hat\nu(a+bn)|^2 \le
C_p^2 \|P_\nu\|_{L^p\to L^q}^2 N^{2(p-1)/p}(2N+1)^{(2-q)/q}.
\end{equation}

(ii) If $2 \le p$, then there exists a constant $C_{q}$ such that 
for every $a \in \mathbb Z$ and $b \in \mathbb N^+$ we have
\begin{equation} \label{sum-estimate2}
\sum_{|n| \le N}|\hat\nu(a+bn)|^2 \le  
C_q^2 \|P_\nu\|_{L^p\to L^q}^2 N^{2/q}(2N+1)^{(p-2)/p}.
\end{equation}
\end{theo}
\begin{proof}
(i) Let $r=q/(q-1)$ be the dual index of $q\le 2$.  For $f \in L^p$ we have $P_\nu f \in L^q$, 
so by the Hausdorff-Young theorem \cite[Theorem XII.2.3]{Z}
\begin{equation} \label{p-q-2}
\Big(\sum |\hat\nu (n)|^r|\hat f(n)|^r\Big)^{1/r} \le
\|P_\nu f\|_q  \le \|P_\nu\|_{L^p \to L^q} \|f\|_p \ .
\end{equation}
Fix $a \in \mathbb Z$ and $b \in \mathbb N^+$ and put $g_N:= \sum_{n=-N}^N \e^{i(a+nb)x} $. Applying \eqref{p-q-2} to $f=g_N$ we obtain
\begin{equation} \label{r-powers}
\Big(\sum_{|n|\le N} |\hat\nu (a+bn)|^r\Big)^{1/r} \le 
	\|P_\nu g_N\|_q  \le \|P_\nu\|_{L^p\to L^q} \|g_N\|_p \ .
\end{equation}
Since $r \ge 2$, we deduce that
\begin{equation} \label{jensen}
\Big(\frac1{(2N+1)}\sum_{|n| \le N} |\hat\nu(a+bn)|^2 \Big)^{1/2} \le
\Big(\frac1{(2N+1)}\sum_{|n| \le N} |\hat\nu(a+bn)|^r \Big)^{1/r} \le
\frac{\|P_\nu\|_{L^p\to L^q}}{(2N+1)^{1/r}} \|g_N\|_p.
\end{equation}

Now $g_N(x)=\e^{iax} D_N(bx)$, so by  change of variable $y=bx$ and periodicity we obtain
$\|g_N\|_p^p =\|D_N\|_p^p$.
When $p>1$, we use the estimate $\|D_N\|_p^p \le C_p N^{p-1} $ \cite[Lemma 2.1]{AAJRS}
in \eqref{jensen} and obtain
$$
\sum_{|n| \le N}|\hat\nu(a+bn)|^2 \le  
\|P_\nu\|_{L^p\to L^q}^2 C_{p}^2 N^{2(p-1)/p}(2N+1)^{(r-2)/r}\  ,
$$
which is \eqref{sum-estimate1}. The computations in \cite{AAJRS} with the estimation 
$\int_0^\infty \big|\frac{\sin x}x\big|^p dx \le 1+\frac{1}{p-1}$  yield the existence of 
an absolute constant $c$ such that $C_p \le \frac{c}{p-1}$ for $1<p<2$.
\smallskip

When $p=1$, \cite[Theorem 16.3.4]{Ed2} yields $\nu<<\mu$ with
$\phi:=\frac{d\nu}{d\mu} \in L^q$. Then $\sum |\hat\nu(n)|^r$ converges, 
by Hausdorff-Young, and using \eqref{jensen} we obtain
$$
\sum_{|n| \le N}|\hat\nu(a+bn)|^2 \le  
\big(\sum_{n \in \mathbb Z} |\hat\nu(n)|^r\big)^{2/r} (2N+1)^{(2-q)/q} \le 
\|\phi\|_q^2 (2N+1)^{(2-q)/q}.
$$
Since $\|P_\nu\|_{L^1\to L^q} =\|\phi\|_q$ by \cite[Exercise VI.9.59]{DS},
 \eqref{sum-estimate1} holds also for $p=1$, with $C_1=1$.
\medskip

(ii) Let $p\ge 2$. By assumption, $ P_\nu^*$ maps  $(L^q)^*$ to $(L^p)^*$. Denoting 
$r=q/(q-1)$ and $s=p/(p-1)$, we have that $P_\nu^*$ maps $L^r$ to $L^s$, with $1< r<s $, 
and $p\ge 2$ implies $s \le 2$.
The dual $P_\nu^*$ is induced by  $\check\nu$, the reflected probability defined by 
$\check\nu(A)=\nu(-A)$. We then have $\widehat{\check\nu}(n)= \overline{\hat\nu(-n)}$, 
and applying \eqref{sum-estimate1} to $P_\nu^*$ yields \eqref{sum-estimate2}, with 
$C_q =C_r \le c(q-1)$.
\end{proof}

{\bf Remarks.} 1. When $P_\nu$ maps $L^p$ to $L^q$ with $1<p<2<q$, we can use
either \eqref{sum-estimate1} with $r=2$ (when $q\ge p/(p-1)$), or \eqref{sum-estimate2} 
with $p=2$ (otherwise). With $t:= \max\{q,p/(p-1)\}$ (observe that $t>2$) we then obtain
$$
\sum_{|n| \le N}|\hat\nu(n)|^2 \le K\cdot N^{2/t}.
$$
Note that in this case, $\limsup_{|n| \to \infty} |\hat\nu(n)|^2 \le \frac2t,$ 
by Hare \cite[Corollaries 2 and 2']{Ha}.
\smallskip

2. When $P_\nu$ maps $L^p$ to $L^q$ with $1<p<q \le 2$, \eqref{r-powers} yields, with
$r=q/(q-1)$,
\begin{equation} \label{r-powers1}
\sum_{|n|\le N} |\hat\nu (n)|^r \le  K N^{(p-1)r/p}.
\end{equation}
Using Abel summation by parts and the estimate 
$\frac1{n^\beta} -\frac1{(n+1)^\beta} < \frac{\beta}{n^{\beta+1}}$, \eqref{r-powers1} yields
 the (necessary) condition \cite[(16.4.9)]{Ed2}. Our proof is different, and seems simpler. 
Moreover,  \eqref{r-powers1} is a better necessary condition, since its failure is enough for 
ruling out mapping $L^p$ to $L^q$.

3. The case $q=2$ of \eqref{sum-estimate1} is in the thesis of Bonami \cite[p. 345]{Bo},
using Fej\'er kernels.

\smallskip

\begin{prop} \label{coefficients} 
Let $\nu$ be a probability on $\mathbb T$. If for some $\alpha>0$
$$
\sum_{n \in \mathbb Z} |\hat\nu(n)|^\alpha < \infty,
$$
then $P_\nu$ is hyperbounded. More precisely, $P_\nu$ maps $L^p$ to $L^2$, for
$p=\max\{1,\frac{2\alpha}{\alpha+2}\}$.
\end{prop}
\begin{proof}
If $\alpha<2$, we have also $ \sum_{n \in \mathbb Z} |\hat\nu(n)|^2< \infty$, and then 
there is $\phi \in L^2$ with $\hat\phi=\hat\nu$, so by uniqueness $\nu<<\mu$ with
$\frac{d\nu}{d\mu}=\phi$, and by Proposition \ref{Lr-deriv} $P_\nu$ maps $L^1$
to $L^2$, so it is hyperbounded.

Assume now $\alpha>2$, and put $p=\frac{2\alpha}{\alpha+2}$. Then $1<p<2$. 
Let $f \in L^p$. By the Hausdorff-Young theorem we have 
$\sum_{n\in\mathbb Z} |\hat f(n)|^q < \infty$, with 
$q=\frac p{p-1}=\frac{2\alpha}{\alpha -2}$ the dual index to $p$. 
We use H\"older's inequality with $s=\frac q2>1$ and 
$\frac1r=1 -\frac1s =1-\frac2q=\frac2{\alpha}$ and obtain
$$
\sum_{n\in\mathbb Z}|\widehat{P_\nu f}(n)|^2=
\sum_{n\in\mathbb Z}|\hat\nu(n)\hat f(n)|^2 \le 
\Big( \sum_{n\in\mathbb Z}|\hat\nu(n)|^{2r}\Big)^{1/r}
\Big( \sum_{n\in\mathbb Z}|\hat f(n)|^{2s}\Big)^{1/s} < \infty,
$$
since $2s=q$ and $2r=\alpha$. We conclude that $P_\nu$ maps $L^p$ to $L^2$, so $P_\nu$ is
hyperbounded.
\end{proof}

{\bf Remarks.} 1. Proposition \ref{coefficients} was first proved by Hare \cite[Corollary 1]{Ha}.
Our proof is different, and seems simpler.

2. Proposition \ref{coefficients} applies whenever 
$|\hat\nu(n)|= \mathcal O(|n|^{-c})$ for some $c>0$. This condition was shown to
imply hyperboundedness of $P_\nu$ in \cite[Theorem XII.5.24]{Z}; 
for additional information see \cite[Theorem 16.4.6(3)]{Ed2}.
 Singular probabilities satisfying this latter condition were constructed by Littlewood and by 
Wiener and Wintner; see \cite[Theorem XII.10.12]{Z}, and the historical section in \cite{BH}.

3. If $\nu$ satisfies the assumptions of Proposition \ref{coefficients}, then $P_\nu$ 
is Harris recurrent, since for $k> \alpha$ we have 
$\sum_n|\widehat{\nu^k}(n)|=\sum_n|\hat\nu(n)|^k<\infty$, so $\nu^k << \mu$, and 
$P_\nu^k=P_{\nu^k}$ is an integral operator.

4. Let $1<p_0 <2$. A sufficient condition of Edwards \cite[p. 303]{Ed2} for $P_\nu$ to
map $L^{p_0}$ to $L^2$ is $|\hat\nu(n)| = \mathcal O(|n|^{-1/s})$ for $|n|>1$, where
$\frac1s=\frac1{p_0} -\frac12$. For this estimate, Proposition \ref{coefficients} yields 
only that $P_\nu$ maps $L^p$ to $L^2$ for $p_0<p<2$ (by taking $\alpha>s$ close to $s$).
\medskip

{\bf Notations.}  For a sequence of numbers $(a_n)_{n\ge N}$ we put $\Delta a_n:=a_n-a_{n+1}$,
and $\Delta^2 a_n:=\Delta(\Delta a_n)=a_n+a_{n+2}- 2a_{n+1}$.

We also denote the Dirichlet Kernel 
$D_n(x)=\frac12+\sum_{k=1}^n\cos(kx)=\frac12\sum_{0\le |k|\le n} \e^{ikx}$ and the Fej\'er
kernel $K_n(x)=\frac1{n+1}\sum_{k=0}^n D_k(x)$.

\begin{lem} \label{convex}
Let $(a_n)_{n\ge N}$ be a positive sequence with $a_n\to 0$, and $\Delta^2 a_n\ge 0$ for $n\ge N$;
 then there exists a constant $C=C((a_n), N)>0$, such that $ C+\sum_{n=N}^\infty a_n\cos(nx)$
is the Fourier series of a non-negative function  function in $ L^1(-\pi,\pi)$.
\end{lem}
\begin{proof}
Given  a real sequence $(a_n)_{n\ge 0}$, consider the formal Fourier series 
$s(x)=\frac12{a_0}+\sum_{n=1}^\infty a_n\cos(nx)$. By twice summation by parts we obtain that, 
formally, when $a_n\to 0$, 
$$
s(x)=\sum_{n=0}^\infty (n+1)\Delta^2 a_n K_n(x)
$$
(see \cite[formula V.(1.7), p.183]{Z}; convexity is not used, only $|D_n(x)| \le \pi/|x|$ 
for $|x|>0$ \cite[Section II.5]{Z}, and $0 \le nK_n(x) \le 1/2\sin^2(\frac12 x)$ 
\cite[formula III.(3.2)]{Z}).

If we assume, in addition to $a_n\to 0$, that $(a_n)_{n\ge 0}$ is convex, i.e. 
$\Delta^2a_n \ge 0$ for $n \ge 0$, then $a_n \ge 0$ for $n\ge 0$ \cite[Theorem III.4.1]{Z}, 
$s(x)$ converges for every  $x \ne 0$ to a non-negative function (still denoted $s$), and  
$\sum_{n=0}^\infty (n+1) \Delta^2 a_n<\infty$; hence $s(x)\in L^1(-\pi,\pi)$ 
\cite[Theorem V.1.5]{Z}. 
\smallskip


Now, under the assumptions of the lemma,  we define $a_0,a_1,\dots,a_{N-1}$ by  
$a_{N-1}=2a_N - a_{N+1}$, $a_{N-2}=2a_{N-1}- a_N$, etc., so that  $\Delta^2 a_n \ge 0$
for $n \ge 0$. Positivity of $s(x)$ yields that 
$\sum_{n=N}^\infty a_n\cos(nx) + \sum_{n=0}^{N-1}a_n\cos(nx)$ is positive for $x \ne 0$, 
so adding the norm 
$$C :=\sup_{|x| \le \pi} \Big|\sum_{n=0}^{N-1}a_n \cos(nx)\Big|,$$  
the series  $C+\sum_{n=N}^\infty a_n\cos(nx)=C+\frac12 \sum_{|n| \ge N} a_n\e^{inx}$ 
is positive for $x\ne 0$. Integrability clearly holds.
\end{proof}

\begin{prop} \label{product=hyper}
There exist two probabilities $\nu_1$ and $\nu_2$ on $\mathbb T$ such that 
$P_{\nu_1}P_{\nu_2}$ is hyperbounded, but both $P_{\nu_1}$ and $P_{\nu_2}$ 
are not hyperbounded.
\end{prop}
\begin{proof}
By Lemma \ref{convex}, the series $g(x)=\sum_{n=2}^\infty \frac{\cos(nx)}{\log n}$ 
and $f(x)=\sum_{n=1}^\infty \frac{\cos(nx)}{\log(2n)}$, are integrable, and
for some positive constant $C_2>0$, $f_2(x)=f(2x)+C_2$ is positive integrable.
Set 
$$h(x)=g(x) - f(2x) =  \sum_{k=1}^\infty \frac{\cos((2k+1)x)}{\log(2k+1)}.$$
Then $h$ is integrable. We would like to show that for some positive constant $C_1$, $h+C_1$ 
is non-negative.

We use Theorem V.2.17 of \cite[p. 189]{Z}. It yields directly that 
$g(x)\approx \frac{\pi}{2x\log^2(1/x)}$ for $x\to 0^+$.  We then apply it with 
$b(u)=\frac1{\log (2u)}$,
and obtain that  $f(2x)\approx\frac{\pi}{4x\log^2(1/x)}$ for $x\to 0^+$. 
Thus $h(x)\approx \frac{\pi}{2x\log^2(1/x)}$ for $x\to 0^+$, so $h(x)>0$ for
$x \in(0,\epsilon)$.  By \cite[Theorem I.2.6]{Z}, the series $g(x)$ and $f(x)$
converge uniformly for $|x|\ge \epsilon$; hence $h(x)$ is bounded below on $|x| \ge \epsilon$.
Finally we conclude that for some positive constant $C_1$, $h+C_1$ is non-negative integrable. 

Let $f_1(x):= h(x)+C_1$, and define $\nu_1$ and $\nu_2$ by 
$\frac{d\nu_j}{d\mu} =\frac{f_j}{\|f_j\|_1}$. Since $\hat\nu_1(n)=0$ for even $n \ne 0$
and $\hat\nu_2(n)=0$ for odd $n$, we have $\nu_1*\nu_2=\mu$, so $P_{\nu_1}P_{\nu_2}$
is hyperbounded. However, as in the proof of  Proposition \ref{not-hyper2}, for any 
$1<p<2$, each $P_{\nu_j}$ does not satisfy  \eqref{hr}, so is not hyperbounded.
\end{proof}

{\bf Remark.} Proposition \ref{product=hyper} shows that powers in Proposition \ref{ritter} 
cannot be  replaced by products.

\begin{prop} \label{monotone}
Let $\nu$ be a probability on $\mathbb T$, such that for some $C>0$,
\begin{equation} \label{korner}
|\hat\nu(n)|^2 \le C \cdot \frac1n \sum_{k=n+1}^{2n} |\hat\nu(k)|^2  \qquad n=1,2,\dots
\end{equation}
If $P_\nu$ is hyperbounded, then $P_\nu$ is Harris recurrent.
More precisely, some convolution power $\nu^k$ is absolutely continuous.
\end{prop} 
\begin{proof} Since $P_\nu$ is hyperbounded, $\nu$ is continuous by Proposition 
\ref{continuous}. For any continuous $\nu$, Wiener's characterization of continuity 
with \eqref{korner} imply that $\nu$ is Rajchman.

Let $P_\nu$ map $L^p$ into $L^q$, $1 \le p <q$.
We use condition \eqref{korner} and the estimate of the Fourier-Stieltjes coefficients 
given by Theorem \ref{sum-estimates}.  

When $1\le p<q\le 2$, \eqref{sum-estimate1} yields
\begin{equation} \label{lower}
N|\hat\nu(N)|^2 \le C^2\sum_{n=N+1}^{2N} |\hat\nu(n)|^2 \le 
C^2\sum_{n=1}^{2N} |\hat\nu(n)|^2 \le K^2(\log N)^{2/p}N^{2(p-1)/p} N^{(r-2)/r},
\end{equation}
where $r =q/(q-1)$. This  yields 
\begin{equation} \label{smallest}
|\hat\nu(N)|^2 \le K^2 (\log N)^{2/p} N^{^{1-\frac2p + 1 -\frac2r }}.
\end{equation}
Since $p<q \le 2$, we have $s:= \frac2p +\frac2r - 2=\frac2p +\frac{2(q-1)}q -2 =
\frac2p- \frac2q >0$.  Hence for $k =[\frac1s]+1 $ we obtain
$$
\sum_{ N=1}^\infty |\widehat{\nu^k}(N)|^2 =
\sum_{ N=1}^\infty |\hat{\nu}(N)|^{2k} < \infty,
$$
so $\nu^k<<\mu$, with $\frac{d\nu}{d\mu} \in L^2$.
In particular, if $P_\nu$ maps $L^1$ into $L^2$, then $\nu*\nu$ is absolutely continuous.
\smallskip

When $2 \le p<q$, \eqref{sum-estimate2} yields similarly
$$
N|\hat\nu(N)|^2 \le C^2 \sum_{n=1}^{2N} |\hat\nu(n)|^2 \le 
K^2(\log N)^{2/r}N^{2(r-1)/r} N^{(p-2)/p},
$$
which yields
$$
|\hat\nu(N)|^2 \le K^2 (\log N)^{2/r} N^{^{1-\frac2p + 1 -\frac2r }}.
$$
As before, $\sum_{ N=1}^\infty |\widehat{\nu^k}(N)|^2 < \infty$ for $k$ and $s$
defined as above, so $\nu^k <<\mu$.
\end{proof}

{\bf Remarks.} 1.  By Cauchy-Schwarz, if 
$|\hat\nu(n)| \le C\cdot \frac1n \sum_{k=n+1}^{2n} |\hat\nu(k)|$ 
for $n\ge 1$, then \eqref{korner} holds.

2. Assumption \eqref{korner} is inspired by \cite{Ko}. 
The averaging is done in order to allow some coefficients to be zero.

3. With a minor change in the proof, assumption \eqref{korner} may be replaced by
\begin{equation} \label{korner2}
|\hat\nu(n)|^2 \le C \cdot \frac1n \sum_{k=1}^{n} |\hat\nu(k)|^2  \qquad n=1,2,\dots\ ,
\end{equation}
which includes the case of $(|\hat\nu(n)|)_{n\ge 1}$ non-increasing. If
$|\hat\nu(n)| \le C \cdot \frac1n \sum_{k=1}^{n} |\hat\nu(k)|$ for $n \ge 1$, then
\eqref{korner2} holds.
%

\bigskip

{\it Application to uniform distribution modulo 1}

Recall \cite[p. 1]{KN} that a sequence of reals $(x_k)_{k\ge 1}$ is called {\it uniformly 
distributed (u.d.) modulo 1} if for any subinterval $0\le a <b \le 1$ its fractional parts 
$\{x_k\}:=x_k-[x_k]$ satisfy
$$
\lim_{N\to \infty} \frac1N\sum_{k=1}^N1_{[a,b)}(\{x_k\}) = b-a.
$$
Since we are concerned with $[0,1)$ (and not an interval of length $2\pi$), we denote now
$e(x):= \e^{2\pi ix}$. By change of variables, now $\hat\nu(n) = \int_0^1 e(-nx)d\nu(x)$
for a probability $\nu$ on $\mathbb T$.

By Weyl's criterion \cite[p. 7]{KN}, $(x_k)_{k \ge 1}$ is u.d. if and only if for every 
integer $m\not=0$
$$
\frac1N\sum_{k=1}^N e(m x_k)\to 0.
$$
Weyl proved that if $(n_k)_{k\ge 1}$ is a sequence of distinct integers, then for Lebesgue
almost every $x$ the sequence $(n_k x)$ is u.d. mod 1 \cite[p. 32]{KN}.

\begin{prop} \label{lyons}
Let $\nu$ be a probability on $\mathbb T$. If 
$|\hat\nu(n)| =\mathcal O((\log|n|)^{-\gamma})$ for some $\gamma>1$, 
then  for every strictly increasing sequence of positive integers 
$(n_k)_{k\ge 1}$, the sequence $(n_kx)$ is u.d. for $\nu$ a.e. $x$. 
\end{prop}
This is a special case (with $\Phi(n)=(\log n)^{-\gamma}$) of Theorem 3 of Lyons \cite{Ly}.
\smallskip




\begin{proof} Fix a probability $\nu$ and a sequence sequence $(n_k)_{k\ge 1}$ of 
distinct integers.  In order to prove that $(n_k x)$ is u.d. mod 1 for $\nu$-a.e. $x$, 
we have to prove that $\nu$-a.e. $x$ satisfies
$$
S_N(x,m):=\frac1N\sum_{k=1}^N e(m n_kx)\to 0 \quad 
\text{for every } \ 0\ne m \in \mathbb Z
$$

Applying to $\nu$ the method of Davenport, Erd\"os, LeVeque  \cite{DEL} 
(see \cite[p. 33]{KN}), it suffices to check that for every $0 \ne m \in \mathbb Z$, 
\begin{equation} \label{del}
\sum_{N=1}^\infty \frac1{N} \int |S_N(x,m)|^2d\nu(x) =
\sum_{N=1}^\infty \frac1{N^3} \sum_{k,j=1}^N \hat\nu(m(n_k-n_j))<\infty.
\end{equation}
When $|\hat\nu(n)|=\mathcal O((\log|n|)^{-\gamma})$, convergence in \eqref{del} holds.
\end{proof}

{\bf Remarks.} 1. If $\nu$ satisfies all the assumptions of Proposition \ref{monotone}, 
then by \eqref{smallest} $|\hat\nu(n)| = \mathcal O(|n|^{-c})$ for some $c>0$.
so we can apply Proposition \ref{lyons}.

2. Proposition \ref{lyons} applies to the above mentioned examples by Littlewood and by
Wiener-Wintner, of $\nu$ singular with $|\hat\nu(n)|=\mathcal O(|n|^{-c})$ for some $c>0$.

3, Additional singular probabilities which satisfy the assumption of Proposition 
\ref{lyons} can be obtained from Theorem 1.2 of K\"orner \cite{Ko} 
(with $\phi(n) =(\log(n+1))^{-\gamma}$).

4. Let $\nu_0$ be the  probability of Proposition \ref{not-hyper2}, and $\nu:=\nu_0*\nu_0$.
Then Proposition \ref{lyons} applies, but $P_\nu$ is not hyperbounded, by Proposition
\ref{ritter}.

\begin{prop} \label{bded-gaps}
Let $\nu$ be a probability on $\mathbb T$, such that $P_\nu$ is hyperbounded.
Then for every  sequence of distinct positive integers $(n_k)_{k\ge 1}$
with bounded gaps, $(n_kx)$ is u.d. for $\nu$ a.e. $x$. 
\end{prop}
\begin{proof} 
Let $(n_k)$ have  bounded gaps: $|n_{k+1}-n_{k}| \le d$ for any $k \ge 1$. 

 By the Cauchy-Schwarz inequality
\begin{equation} \label{double}
\sum_{k,j=1}^N\hat\nu(m(n_k-n_j)) \le
N|\hat\nu(0)|+ N \Big( \sum_{k\not=j}^N |\hat\nu(m(n_k-n_j))|^2 \Big)^{1/2}.
\end{equation}

{\it Claim: Each value of $n_k-n_j$ can not occur more than $N$ times.}

\noindent
{\it Proof of claim:} Denote $c_N(\ell) := |\{1\le k \le N: n_k=\ell\}|$.
Since $(n_k)$ are distinct, $c_N(\ell)$ is 0 or 1, and is 1 for $N$ values of $\ell$. Put 
$n_N^*= \max_{1\le k \le N} n_k$ and $V_N(t):=|\{1 \le k<j \le N: n_k-n_j=t\}|$. Then
$$
V_N(t)=\sum_{\ell=1}^{n_N^*} c_N(\ell+t)c_N(\ell) \le N.
$$
This means that the value $t \ne 0$ occurs as a difference at most $N$ times.
\smallskip

The indices of the Fourier-Stieltjes coefficients in \eqref{double} are included
in the interval with end points the maximal difference, which is at most $|m|(N-1)d$ 
(note that $m$ is fixed), so
\begin{equation} \label{double-sum}
\sum_{k,j=1}^N\hat\nu(m(n_k-n_j)) \le
N|\hat\nu(0)|+N \Big(N \sum_{k=-|m|Nd}^{|m|Nd} |\hat\nu(k)|^2 \Big)^{1/2} .
\end{equation}

Since $P_\nu$ is hyperbounded, we apply Theorem \ref{sum-estimates}.
We first assume that $P_\nu$ maps $L^p$ into $L^q$ with $1 \le p<2$, and $q \le 2$. Then
\eqref{double-sum} yields
$$
\sum_{k,j=1}^N\hat\nu(m(n_k-n_j)) \le
CN +C'_m(\log N)^{1/p}N^{3/2}N^{(p-1)/p}N^{(r-2)/2r},
$$
by the estimate \eqref{sum-estimate1}, with $r=q/(q-1)$.
But 
$$
\frac32 +\frac{p-1}p +\frac{r-2}{2r}=3 -\frac1p -\frac{q-1}q =2 -\frac1p+\frac1q <2
$$
since $\frac1p>\frac1q$, so substituting into \eqref{del} we obtain convergence of
the series; hence $(n_kx)$ is uniformly distributed mod 1 for $\nu$-a.e. $x$.

We now assume that $P_\nu$ maps $L^p$ to $L^q$ with $2 \le p<q$. Then 
\eqref{double-sum} yields
$$
\sum_{k,j=1}^N\hat\nu(m(n_k-n_j)) \le
CN +C'_m(\log N)^{1/r   }N^{3/2}N^{(r-1)/r}N^{(p-2)/2p},
$$
by the estimate \eqref{sum-estimate2}, with $r=q/(q-1)$.
But 
$$
\frac32 +\frac{r-1}r +\frac{p-2}{2p}=3 -\frac1r -\frac1p =2 -\frac1p+\frac1q <2,
$$
and we obtain as before convergence of the series \eqref{del}, which proves that
 $(n_kx)$ is uniformly distributed mod 1 for $\nu$-a.e. $x$.
\end{proof}
 
{\bf Remark.} In the right-hand side of the estimate \eqref{double-sum}, we can replace
  the endpoints of the summation interval by $\pm |m|n_N^*$. 
When $P_\nu$ maps $L^p$ into $L^q$, there is an $\alpha= \alpha(p,q)$
such that $(n_k x)$ is u.d. for $\nu$ a.e. $x$ when  $(n_k)$ are distinct with 
$n_N^*=\mathcal O(N^\alpha)$; in particular, this applies when $n_k=\mathcal O(k^\alpha)$.
When $1 \le p<q \le 2$ or $2\le p<q$, Theorem \ref{sum-estimates} yields
$\alpha(1-\frac2p+\frac2q) <1$ in either case. When $p=1$ and $q=2$, or
$p=2$ and $q=\infty$, for any $\alpha \ge 1$  the sequence $(n_k x)$ is u.d. mod 1.
Since $\frac1p > \frac1q$, in all other cases $1-\frac2p+\frac2q <1$, so
the case of bounded gaps (for which $\alpha=1$) is always a special case.

\medskip

\bigskip

{\bf Acknowledgement.} The second author expresses his deep gratitude to the late Professor 
Shaul Foguel (1931-2020), his Ph.D. advisor, for introducing him to research in the ergodic 
theory of Markov operators.

\end{document}